\documentclass[twoside]{amsart}
\copyrightinfo{}{}

\usepackage{amssymb}
\usepackage{amsmath}
\usepackage{cite}
\usepackage{mathrsfs}

\makeatletter

\address{\bigskip\hfil\begin{tabular}{l@{}}
            School of Mathematics and Statistics F07\\  
            University of Sydney,
            Sydney N.S.W. 2006.\hfill\qquad
                        {\tt mathas@maths.usyd.edu.au}\\  
            Australia.\hfill\qquad
                    {\tt www.maths.usyd.edu.au/u/mathas/}
          \end{tabular}}

\let\atop\@@atop

\def\And{\text{\ and\ }}

\def\For{\text{\ for\ }}

\def\ForSome{\text{\ for some\ }}
\def\If{\text{\ if\ }}

\def\Otherwise{\text{\ otherwise}}
\def\th{{\text{th}}}

\def\Set[#1]#2|#3|{\Big\{\ #2\ \Big| \
            \vcenter{\hsize #1mm\centering#3}\Big\}}

{\catcode`\|=\active
  \gdef\set#1{\mathinner{\lbrace\,{\mathcode`\|"8000%
                                   \let|\midvert #1}\,\rbrace}}
}
\def\midvert{\egroup\mid\bgroup}

\def\Number#1{\refstepcounter{equation}
              \leqno(\theequation)\if*#1%
              \else\def\@currentlabel{{\rm\theequation}}\label{#1}%
              \fi}
\def\Dag{\ifmmode\leqno(\dag)\else$(\dag\)$\fi}
\def\DDag{\ifmmode\leqno(\ddag)\else$(\ddag\)$\fi}

\numberwithin{equation}{section}

\swapnumbers
\newtheorem{Definition}[equation]{Definition}
\newtheorem{Theorem}[equation]{Theorem}
\newtheorem{Proposition}[equation]{Proposition}
\newtheorem{Lemma}[equation]{Lemma}
\newtheorem{Corollary}[equation]{Corollary}

\theoremstyle{remark}
\newtheorem{Remark}[equation]{Remark}
\newtheorem{Remarks}[equation]{Remarks}
\newtheorem{Example}[equation]{Example}

\newenvironment{Point}[2]%
  {\ifx*#2\let\pointlabel\relax\else\def\pointlabel{#2}\fi
   \refstepcounter{equation}\trivlist
   \item[\hskip\labelsep\bf\theequation
         \ifx\pointlabel\relax\else\space\pointlabel\space\fi]
   \ignorespaces #1
  }{\relax}

\let\DS\displaystyle
\let\MC\multicolumn
\let\SR\stackrel
\def\Prod{\displaystyle\prod}
\def\Sum{\displaystyle\sum}

\def\cf{{\it cf.\space}}

\def\){\big)}
\def\({\big(}
\let\>\rangle
\let\<\langle

\let\incl\hookrightarrow
\let\ss\subseteq
\let\0\varnothing              

\let\realb@r\bar
\let\bar\overline

\let\gedom\trianglerighteq
\let\gdom\vartriangleright

\def\Z{{\mathbb Z}}
\def\Q{{\mathbb Q}}
\def\R{{\mathbb R}}
\def\C{{\mathbb C}}

\let\To\longrightarrow
\def\map#1#2{\,{:}\,#1\!\longrightarrow\!#2}
\def\mapsto{\!\longmapsto\!}

\title[Matrix units and generic degrees]%
      {Matrix units and generic degrees for the Ariki--Koike algebras}
\author{Andrew Mathas}


\makeatletter
\let\bar\realb@r
\makeatother

\def\len(#1){\ell(#1)}
\def\res{\mathop{\rm res}\nolimits}
\def\Sym{\mathfrak S}

\def\Q{\mathbf q}

\def\c{{\mathbf c}}

\def\s{\mathfrak s}
\def\t{\mathfrak t}
\def\u{\mathfrak u}
\def\v{\mathfrak v}
\def\a{\mathfrak a}
\def\b{\mathfrak b}

\def\defeq{\overset{\text{def}}=}

\def\R{\mathcal R}
\def\ZZ{\mathcal Z}
\def\H{\mathscr H}
\def\HZ{\H_\ZZ}
\def\Hlam{\check\H^\lambda}

\def\PS{^{\phantom*}}

\def\barl{{\bar\lambda}}
\def\wbarl{{w_{\bar\lambda}}}
\def\lamp{{\lambda'}}
\def\tlam{{\t^\lambda}}
\def\tllam{{\t_\lambda}}

\def\tlamp{{\t^\lamp}}
\def\Flam{F_\tlam}

\def\glam{\gamma_\tlam}

\def\ulam{u_\lambda}

\def\standard{\text{\ standard}}
\def\rest{{\downarrow}}

\def\Std{\operatorname{\rm Std}}
\def\Shape{\operatorname{\rm Shape}}
\def\rtuple#1{(#1^{(1)},\dots,#1^{(r)})}

\usepackage[enableskew,vcentermath]{youngtab}
\def\bitab(#1|#2){\bigg(\ \young(#1\ ,\quad\young(#2)\ \bigg)}
\def\tritab(#1|#2|#3){{\bigg(\begin{array}[b]{*3c}\young(#1)\ ,&\young(#2)\ ,%
                                    &\young(#3)\end{array}\bigg)}}

\begin{document}
\begin{abstract}

We compute the generic degrees of the Ariki--Koike algebras by first
constructing a basis of matrix units in the semisimple case. As a
consequence, we also obtain an explicit isomorphism from any
semisimple Ariki--Koike algebra to the group algebra of the
corresponding complex reflection group.
\end{abstract}
\maketitle

\section{Introduction}

The cyclotomic Hecke algebras were introduced by  Ariki and
Koike~\cite{AK1,Ariki:hecke} and Brou\'e and Malle~\cite{BM:cyc}. It
is conjectured~\cite{BM:cyc} that these algebras play a r\^ole in the
representation theory of reductive groups similar to (but more
complicated than) that played by the Iwahori--Hecke algebras (see, for
example, \cite{Carter:2}).  In particular, it should be possible to
use these algebras to compute the degrees (and more generally
characters) of certain representations of reductive groups; more
precisely, we can attach a polynomial to each irreducible
representation~$\H$, called its generic degree, and appropriate
specialisations of this polynomial should compute the dimensions of
the corresponding irreducible representations of the finite groups of
Lie type.

The purpose of this paper is to compute the generic degrees for the
cyclotomic Hecke algebras of type $G(r,1,n)$; these polynomials have
also been computed by Geck, Iancu and Malle~\cite{GIM}.  Further
results of Malle~\cite{Malle:gendeg} and Malle and the
author~\cite{MM:trace} give the generic degrees for all of the
cyclotomic Hecke algebras corresponding to imprimitive complex
reflection groups.  Malle~\cite{Malle:expgendeg} has recently computed
the generic degrees for cyclotomic algebras for the primitive complex
reflection groups (modulo the assumption that the corresponding
cyclotomic Hecke algebras are symmetric), so this completes the
calculation of the generic degrees of the cyclotomic Hecke algebras
associated with complex reflection groups. 

Two important special cases of the Ariki--Koike algebras are the
Iwahori--Hecke algebras of types $A$ and $B$; the generic degrees of
these algebras are well--known and were first computed by
Hoefsmit~\cite{Hoefsmit}. Later Murphy~\cite{murphy:hecke} gave an
easier derivation of Hoefsmit's formulae for the generic degrees of
the Iwahori--Hecke algebras of type~$A$ using different, but related,
techniques. 

This article is largely inspired by Murphy's
paper~\cite{murphy:hecke}; however, with hindsight we are able to take
quite a few shortcuts.  Along the way we give a new and quite elegant
treatment of the representation theory of the semisimple Ariki--Koike
algebras.  In particular, we explicitly construct the primitive
idempotents and the matrix units in the Wedderburn decomposition of
$\H$. One of the nice features of our approach is that we use the
modular theory (more accurately, the cellular theory) to understand
the semisimple case.  Another advantage of our approach is that our
calculation of the generic degrees is almost entirely algebraic and it
is only near of the end of the paper that we need to introduce
combinatorial arguments (however, this is not to say that our approach
is painless; in addition, showing that our result agrees with
\cite{GIM} is a long and unpleasant calculation).

\section{The Ariki--Koike algebras}

Let $R$ be a commutative domain with $1$ and fix elements $q$,
$Q_1,\dots,Q_r$ in $R$ with~$q$ invertible. Let
$\Q=(q;Q_1,\dots,Q_r)$. The {\sf Ariki--Koike} 
algebra $\H=\H_\Q(n)$ is the unital associative
algebra with generators $T_0,T_1,\dots,T_{n-1}$ and relations 
$$\begin{array}{rlll} (T_0-Q_1)\cdots(T_0-Q_r) &=&0, \\
  T_0T_1T_0T_1&=&T_1T_0T_1T_0,\\ (T_i+q)(T_i-1) &=&0,&
            \text{for $1\le i\le n-1$,} \\ 
  T_{i+1}T_iT_{i+1}&=&T_iT_{i+1}T_i,&
            \text{for $1\le i\le n-2$,}\\
  T_iT_j&=&T_jT_i,& \text{for $0\le i<j-1\le n-2$.} 
\end{array}$$
The Ariki--Koike algebra is a deformation of the group algebra of
$W_{r,n}$, where $W_{r,n}=(\Z/r\Z)\wr\Sym_n$ is the wreath product of the
cyclic group of order $r$ and the symmetric group $\Sym_n$ of degree $n$.

Define elements $L_m=q^{1-m}T_{m-1}\ldots T_1T_0T_1\ldots T_{m-1}$ for
$m=1,2,\ldots,n$; these are analogues of the $q$--Murphy operators of
the Iwahori--Hecke algebras of type~$\bf A$
\cite{DJ:blocks,murphy:hecke}. An easy calculation using the relations
in $\H$ (see \cite[3.3]{AK1} and \cite[(2.1), (2.2)]{DJ:blocks}) shows
that we have the following results.

\begin{Point}
{\it}* Suppose that $1\le i\le n-1$ and $1\le m\le n$. Then
\begin{enumerate}
\item $L_i$ and $L_m$ commute.
\item $T_i$ and $L_m$ commute if $i\ne m-1, m$.
\item $T_i$ commutes with $L_iL_{i+1}$ and $L_i+L_{i+1}$.
\item If $a\in R$ and $i\ne m$ then $T_i$ commutes with 
$(L_1-a)(L_2-a)\ldots(L_m-a)$.
\end{enumerate}\label{L-comm}\end{Point}

Using the elements $T_w$ and $L_m$ defined above, Ariki and Koike
proved the following result which gives a basis for~$\H$.

\begin{Point}
{\it}{(Ariki--Koike \cite[(3.10)]{AK1})}
The algebra $\H$ is free as an $R$--module with basis
$$\set{L_1^{c_1}L_2^{c_2}\ldots L_n^{c_n}T_w|%
       w\in\Sym_n \And 0\le c_m\le r-1\For m=1,2,\ldots,n}.$$
In particular, $\H$ is free of rank $r^nn!$
\label{AK basis}
\end{Point}

We call this basis the {\sf Ariki--Koike basis} of $\H$. Let $\H(\Sym_n)$
be the subalgebra of $\H$ generated by $T_1,\dots,T_{n-1}$. It follows
from (\ref{AK basis}) that $\H(\Sym_n)$ is isomorphic to the
Iwahori--Hecke algebra of $\Sym_n$ and is free as an $R$--module with
basis $\set{T_w|w\in\Sym_n}$.

Let $\tau\map\H R$ be the $R$--linear map determined by 
$$\tau(L_1^{c_1}L_2^{c_2}\ldots L_n^{c_n}T_w) 
           =\begin{cases} 1,&\If c_1=\dots=c_n=0\And w=1,\\
	           0,&\Otherwise,
           \end{cases}$$
where $w\in\Sym_n$ and $0\le c_i<r$ for $i=1,\dots,n$.

The function $\tau$ was introduced by Bremke and
Malle~\cite{BM:redwds} who showed that $\tau$ is a trace form  and
that $\tau$ is essentially independent on the choice of basis of
$\H$. It is not obvious from the definition above that $\tau$
coincides with the form introduced by Bremke and Malle; however,
this was proved by Malle and the author in~\cite{MM:trace} where we
also showed that  $\tau$ is non--degenerate whenever $Q_1,\dots,Q_r$
are invertible in $R$. 

For future reference we note the following two important properties
of~$\tau$; the first is Bremke and Malle's result that $\tau$ is a
trace form and the second follows easily from the definition and
well--known properties of the trace form $\tau$ in the case $r=1$
(see, for example, \cite[Prop.~1.16]{M:ULect}).

\begin{Point}
{\it}*
\begin{enumerate}
\item Suppose $h_1,h_2\in\H$. Then $\tau(h_1h_2)=(h_2h_1)$.
\item Suppose that $x,y\in\Sym_n$ and that $0\le c_i<r$ for 
$i=1,\dots,n$. Then
$$\DS\tau(L_1^{c_1}L_2^{c_2}\ldots L_n^{c_n}T_xT_y)
      =\begin{cases} q^{\len(x)},&\If c_1=\dots=c_n=0\And x=y^{-1},\\
                        0,&\Otherwise.
\end{cases}$$
\end{enumerate}\label{tau properties}\end{Point}

In this paper we will mainly be concerned with the semisimple Ariki--Koike
algebras; these were classified by Ariki~\cite{Ariki:ss} who showed
that when $R$ is a field $\H$ is semisimple if and only if
$$P_\H(\Q)=\prod_{i=1}^n (1+q+\dots+q^{i-1})\cdot
  \prod_{1\le i<j\le r}\prod_{-n<d<n}(q^dQ_i-Q_j)\Number{ppoly}$$
is a non--zero element of $R$. For most of what we do it will be
enough to assume that~$R$ is a ring in which $P_\H(\Q)$ is invertible.

A {\sf multipartition} of $n$ is a ordered $r$--tuple
$\lambda=\rtuple\lambda$ of partitions $\lambda^{(s)}$ such that
$n=\sum_{s=1}^r|\lambda^{(s)}|$; we write $\lambda\vdash n$. In the
semisimple case Ariki and Koike~\cite{AK1} constructed an irreducible
$\H$--module $S^\lambda$, called a {\sf Specht module}, for each
multipartition~$\lambda$ of~$n$. Further, they showed that
$\set{S^\lambda|\lambda\vdash n}$ is a complete set of pairwise
non--isomorphic irreducible $\H$--modules.  Let~$\chi^\lambda$ be the
character of $S^\lambda$.

Assume that $\H$ is semisimple. Then $\tau$ can be written as a linear
combination of the irreducible characters of $\H$ because $\tau$ is a
trace form. 

\begin{Definition}
Suppose that $R$ is a field and that $P_\H(\Q)\ne0$. The {\sf Schur
elements} of $\H$ are the elements $s_\lambda(\Q)\in R$ such that 
$$\tau=\sum_\lambda \frac1{s_\lambda(\Q)}\chi^\lambda,$$
where $\lambda$ runs over the multipartitions of $n$.  
\end{Definition}

The rational functions $\frac1{s_\lambda(\Q)}$ are also
called the weights of~$\H$.

The {\sf generic degrees} of $W_{r,n}$ are certain ``spetsial''
specializations of the rational functions $s_\eta(\Q)/s_\lambda(\Q)$,
where $\eta=\((n),(0),\dots,(0)\)$ ---~so $s_\eta(\Q)$ is the Schur
element corresponding to the trivial representation of $\H$
($s_\eta(\Q)=s_{\eta_1}(\Q)$ is computed in Example~\ref{n_t example}). The
spetsial specializations of the rational functions
$s_\eta(\Q)/s_\lambda(\Q)$ are polynomials in $q$ with rational
coefficients; moreover, for these specializations~$s_\eta(\Q)$ is
equal to the Poincar\'e polynomial of the coinvariant algebra of the
reflection representation of $W_{r,n}$.  These results are due to
Malle and can be found in~\cite{Malle:expgendeg,Malle:gendeg}.  

In the special case when $r=1,2$ the group $W_{r,n}$ is a Weyl group
(rather than just a complex reflection group), and here the generic
degrees were first computed by Hoefsmit~\cite{Hoefsmit}; they can be
found, for example, in \cite{Carter:2}.

One of the motivations for writing this paper was to compute the Schur
elements and hence the generic degrees. Ostensibly the Schur elements
depend in a non--uniform way upon the choice of $q,Q_1,\dots,Q_r$;
however, we shall see that in fact they can be expressed as rational
functions in $q,Q_1,\dots,Q_r$ which depend only on~$\lambda$. The
expression we obtain is a generalization of the hook length formula of
Frame, Robinson and Thrall~\cite{FRT} for the dimensions of the
irreducible representations of the symmetric groups.  This is not
unexpected because $\H\cong RW_{r,n}$ when $R=\C$, $q=1$
and~$Q_s=\zeta^s$, so the orthogonality relations for the characters
of~$W_{r,n}$ imply that
$\chi^\lambda(1)=|W_{r,n}|/s_\lambda(1;\zeta^s)$.  As the referee
remarked, it is worth noting that in general~$P_\H(1;\zeta^s)$ is only
a scalar multiple of $|W_{r,n}|=r^nn!$ since $\prod_{1\le i<j\le
r}(\zeta^i-\zeta^j) =(-1)^{\frac12\binom{r-1}2}e^{\frac e2}$ by
\cite[2.22]{Malle:gendeg}.

We will compute the Schur elements of $\H$ by explicitly constructing
a set of primitive idempotents in $\H$ and then applying the following
Lemma (which is really a well--known fact about symmetric algebras).

\begin{Lemma}
Assume that $R$ is a field and that $\H$ is semisimple.
Let $\lambda$ be a multipartition of $n$ and suppose that
$e_\lambda$ is a primitive idempotent in $\H$ such that 
$S^\lambda\cong e_\lambda\H$. Then $s_\lambda(\Q)=\frac
1{\tau(e_\lambda)}$.
\label{generic degree lemma}
\end{Lemma}

\begin{proof}
Suppose first that $R$ is a field of characteristic
zero. Let $E_\lambda$ be the primitive central idempotent
corresponding to the irreducible module $S^\lambda$.  By
definition $\tau=\sum_\mu \frac1{s_\mu(\Q)}\chi^\mu$, where $\mu$ runs
over the multipartitions of~$n$; therefore, 
$$\tau(E_\lambda)=\sum_{\mu\vdash n} \frac1{s_\mu(\Q)}\chi^\mu(E_\lambda)
                 =\frac1{s_\lambda(\Q)}\chi^\lambda(E_\lambda)
                 =\frac1{s_\lambda(\Q)}\chi^\lambda(1).$$
Now $E_\lambda=e_1+\dots+e_N$, where $e_1,\dots,e_N$ are 
primitive idempotents with $e_i\H\cong S^\lambda$, for all $i$,
and $N=\dim S^\lambda=\chi^\lambda(1)$.
The idempotents $e_1,\dots,e_N$ belong to the same Wedderburn
component of $\H$, so there exist invertible elements
$u_i\in\H$ such that $e_i=u_ie_1u_i^{-1}$ for all $i$.  Consequently,
$\tau(e_a)=\tau(u_ie_1u_i^{-1})=\tau(e_1)$ since $\tau$ is a trace
form; hence, $\tau(E_\lambda)=\chi^\lambda(1)\tau(e_1)$. Without
loss,  $e_1=e_\lambda$ so the Lemma follows.

The case where $R$ is a field of positive characteristic now follows
by a specialization argument (using, for example, Theorem~\ref{idempotents})
which we leave to the reader.
\end{proof}

\begin{Example}
Fix $t$ with $1\le t\le r$ and let
$\eta_t=(\eta_t^{(1)},\dots,\eta_t^{(r)})$ be the multipartition
of~$n$ with $\eta_t^{(s)}=(n)$ if $s=t$ and $\eta_t^{(s)}=(0)$
otherwise. We will compute the Schur elements $s_{\eta_t}(\Q)$.  Let
$x_{\eta_t}=\sum_{w\in\Sym_n}T_w$ and $u_{\eta_t}=\prod_{s\ne
t}\prod_{k=1}^n(L_k-Q_s)$ --- in the product, $1\le s\le r$~--- and
set $m_{\eta_t}=u_{\eta_t} x_{\eta_t}=x_{\eta_t}u_{\eta_t}$
(\cf~(\ref{std basis})).  It follows from (\ref{L-comm}) that $u_{\eta_t}$
is central in~$\H$.  Further,  the relations imply that
$T_0u_{\eta_t}=Q_tu_{\eta_t}$ and $T_w
x_{\eta_t}=q^{\len(w)}x_{\eta_t}$ for $w\in\Sym_n$; it follows that
$L_k m_{\eta_t}=q^{k-1}Q_tm_{\eta_t}$ for $k=1,\dots,n$. Therefore,
the module $m_{\eta_t}\H$ is one dimensional and, in particular,
irreducible; in fact, $S^{\eta_t}\cong m_{\eta_t}\H=Rm_{\eta_t}$ (for
example, use (\ref{L_k action})). Moreover, by what we have said
$$m_{\eta_t}^2 =[n]_q^!\prod_{s\ne t}^r
         \prod_{k=1}^n(q^{k-1}Q_t-Q_s)\cdot m_{\eta_t},$$
where $[n]_q^!=\prod_{k=1}^n(1+q+\dots+q^{k-1})$; so $m_{\eta_t}$ is a
scalar multiple of the primitive idempotent which generates $S^{\eta_t}$. 
Hence, by the Lemma,
\begin{align*}
s_{\eta_t}(\Q)
   &=\(\tau(m_{\eta_t})\)^{-1}
     [n]_q^!\prod_{s\ne t}\prod_{k=0}^{n-1}(q^kQ_t-Q_s)\\
   &=(-1)^{n(r-1)}[n]_q^!\prod_{s\ne t}Q_s^{-n}\cdot%
         \prod_{s\ne t}\prod_{k=0}^{n-1}(q^kQ_t-Q_s).
\end{align*}
Similar arguments give the Schur elements for the multipartition
which is conjugate to $\eta_t$; alternatively, they are given by
Corollary~\ref{symmetry} and the calculation above.

There is an action of $\Sym_r$ on the set of multipartitions of $n$
(by permuting components) and also on the rational functions in
$Q_1,\dots,Q_r$ (by permuting parameters). When $\H$ is semisimple the
Specht modules are determined up to isomorphism by the action of
$L_1,\dots,L_n$; as the relation $\prod_{s=1}^r(T_0-Q_s)=0$ is
invariant under the $\Sym_r$--action it follows that
$s_{v\cdot\lambda}(\Q)=v\cdot s_\lambda(\Q)$ for all multipartitions
$\lambda$ and all $v\in\Sym_r$; this is also clear from
Theorem~\ref{idempotents}(i). In the case where $\lambda=\eta_t$ this symmetry
is evident in the formulae above.  \label{n_t example}\end{Example}

\section{An orthogonal basis for $\H$}

If $R$ is a field and $P_\H(\Q)$ is non--zero then $\H$ is a split
semisimple algebra; hence,~$\H$ has a basis which corresponds to the
matrix units in its Wedderburn decomposition. In this section we
explicitly construct a Wedderburn basis for $\H$. We begin by
recalling the standard basis of $\H$ from \cite{DJM:cyc}.

A {\sf multipartition} of $n$ (with $r$ components) is an $r$--tuple
$\lambda=\rtuple\lambda$ of partitions such that
$|\lambda^{(1)}|+\dots+|\lambda^{(r)}|=n$. Let $\Lambda^+$ be the
set of multipartitions of $n$; then $\Lambda^+$ becomes a poset
under {\sf dominance} where $\lambda\gedom\mu$ if for all 
$1\le s\le r$ and all $i\ge1$ we have
$$\sum_{t=1}^{s-1}|\lambda^{(t)}|+\sum_{j=1}^i\lambda^{(s)}_j
       \ge\sum_{t=1}^{s-1}|\mu^{(t)}|+\sum_{j=1}^i\mu^{(s)}_j.$$
We also write $\lambda\gdom\mu$ if $\lambda\gedom\mu$ and
$\lambda\ne\mu$.

The {\sf diagram} of a multipartition $\lambda$ is the set
$$[\lambda]=\set{(i,j,c)|1\le j\le\lambda^{(c)}_i\And1\le c\le r}.$$
We will think of $[\lambda]$ as being the $r$--tuple of diagrams
of the partitions~$\lambda^{(c)}$, for $1\le c\le r$.  A
{\sf $\lambda$--tableau} is a bijection $\t\map{[\lambda]}\{1,2,\dots,n\}$. If
$\t$ is a $\lambda$--tableau write $\Shape(\t)=\lambda$. As with diagrams,
we will think of a tableau $\t$ as an $r$--tuple of tableaux
$\t=\rtuple\t$, where $\t^{(c)}$ is a $\lambda^{(c)}$--tableau.  The
tableaux~$\t^{(c)}$ are called the {\sf components} of $\t$.  A tableau is
{\sf standard} if in each component the entries increase along the rows and
down the columns; let $\Std(\lambda)$ be the set of standard
$\lambda$--tableaux.

We identity a tableau $\t$ with an $r$--tuple of
labelled diagrams; for example
$$\tritab(123,4|56,7|8,9)\qquad\And\qquad
  \tritab(136,2|49,7|5,8) $$
are two $\((3,1),(2,1),(1^2)\)$--tableaux. Both of these tableaux
are standard.

Given a multipartition $\lambda$ let $\tlam$ be the
$\lambda$--tableau with the numbers $1,2,\dots,n$ entered in
order first along the rows of $\t^{\lambda^{(1)}}$ and then the rows
of $\t^{\lambda^{(2)}}$ and so on. For example, the first of the
$\((3,1),(2,1),(1^2)\)$--tableaux above is 
$\t^{((3,1),(2,1),(1^2))}$.

The symmetric group $\Sym_n$ acts from the right on the set of
$\lambda$--tableaux; let
$\Sym_\lambda
     =\Sym_{\lambda^{(1)}}\times\dots\times\Sym_{\lambda^{(r)}}$
be the row stabilizer of $\tlam$. For any $\lambda$--tableau $\t$
let $d(\t)$ be the unique element of $\Sym_n$ such that
$\t=\tlam d(\t)$.

Let $*$ be the $R$--linear antiautomorphism of $\H$ determined by
$T_i^*=T_i$ for all~$i$ with $0\le i\le n-1$. Then $T_w^*=T_{w^{-1}}$ and
$L_k^*=L_k\PS$ for all $w\in\Sym_n$ and for~$k=1,2,\ldots,n$. 

We can now define the cellular basis of $\H$ constructed
in~\cite{DJM:cyc}. Fix a multipartition~$\lambda$ and let 
$\a=(\a_1,\dots,\a_r)$, where $\a_s=|\lambda^{(1)}|+\dots+|\lambda^{(s-1)}|$
for $1\le s\le r$. Define $m_\lambda=x_\lambda\ulam^+$, where
$$ x_\lambda=\sum_{w\in\Sym_\lambda}T_w\qquad\And\qquad
u_\lambda^+=\prod_{s=2}^r\prod_{k=1}^{\a_s}(L_k-Q_s).$$
It follows from (\ref{L-comm}) that $x_\lambda$ and $\ulam^+$
commute.  Finally, given standard $\lambda$--tableaux~$\s$ and
$\t$ let $m_{\s\t}=T_{d(\s)}^*m_\lambda T_{d(\t)}^{\phantom*}.$

Whenever we write $m_{\s\t}$ in what follows $\s$ and $\t$ will be
standard tableaux of the {\it same shape} (and similarly, for
$f_{\s\t}$ etc.).

\begin{Point}
{\it}{(Dipper--James--Mathas~\cite[Theorem 3.26]{DJM:cyc})}
The Ariki--Koike algebra $\H$ is free as an $R$--module with
cellular basis
$\set{m_{\s\t}|\s,\t\in\Std(\lambda)\ForSome\lambda\in\Lambda^+}$.
\label{std basis}\end{Point}

One consequence of this result is that the $R$--module $\Hlam$ which
has as basis the set of $m_{\u\v}$ where
$\Shape(\u)=\Shape(\v)\gdom\lambda$ is a two--sided ideal of $\H$.  It
follows from (\ref{std basis}) that $m_\mu\H\subset\Hlam$ whenever
$\mu\gdom\lambda$.  The {\sf Specht module} $S^\lambda$ is the right
$\H$--module $m_\lambda\H/(m_\lambda\H\cap\Hlam)$, a submodule of
$\H/\Hlam$. Thus, $S^\lambda$ is free as an $R$--module with basis
$\set{m_\t|\t\in\Std(\lambda)}$, where $m_\t=m_{\tlam\t}+\Hlam$.
Further, there is a natural associative bilinear form $\<\ ,\ \>$ on
$S^\lambda$ which is determined by $$\<m_\s,m_\t\>m_\lambda\equiv
m_{\tlam\s}m_{\t\tlam}\mod\Hlam.$$

Before we can begin we need some further notation and a
result from \cite{JM:cyc-Schaper}. If~$\t$ is any
tableau and $k\ge0$ is an integer let $\t\rest k$ be the
subtableau of $\t$ which contains the integers $1,2,\dots,k$.
Observe that $\t$ is standard if and only if $\Shape(\t\rest k)$ is
a multipartition for all~$k$. We extend the
dominance order to the set of standard tableaux by defining
$\s\gedom\t$ if $\Shape(\s\rest k)\gedom\Shape(\t\rest k)$ for
$k=1,2,\dots,n$; again, we write $\s\gdom\t$ if $\s\gedom\t$ and
$\s\ne\t$.

Write $\res_\t(k)=q^{j-i}Q_c$ if $k$ appears in row $i$ and column $j$
of $\t^{(c)}$; then $\res_\t(k)$ is the {\sf residue} of $k$ in $\t$.
We can now state the result we need from \cite{JM:cyc-Schaper}.

\begin{Point}
{\it}{\cite[Prop. 3.7]{JM:cyc-Schaper}}
Let $\s$ and $\t$ be standard $\lambda$--tableaux and suppose
that $k$ is an integer with $1\le k\le n$. Then there exist 
$a_\v\in R$ such that
$$m_{\s\t}L_k=\res_\t(k)m_{\s\t}
       +\sum_{\SR{v\in\Std(\lambda)}{\v\gdom\t}}a_\v m_{\s\v}
           \bmod\Hlam.$$
\label{L_k action}\end{Point}

From our current point of view the importance of this result
derives from the observation that 
$m_{\s\t}\(L_k-\res_\u(k)\)/\(\res_\t(k)-\res_\u(k)\)=m_{\s\t}$ plus
a linear combination of more dominant terms, providing that
$\res_\u(k)\ne\res_\t(k)$; this motivates the
next definition.

Let $\R(k)=\set{q^dQ_s|1\le s\le r,\ |d|<k
                \And d\ne0\If r=1\And k=2,3}$;
be the set of possible residues $\res_\t(k)$ as $\t$ runs
over the standard tableaux. 

\begin{Definition}
[\protect{cf. \cite[Defn. 3.11]{JM:cyc-Schaper}}] 
Suppose that $\s$ and $\t$ are standard $\lambda$--tableaux.
\begin{enumerate}
\item Let 
$F_\t=\Prod_{k=1}^n
  \prod_{\substack{c\in\R(k)\\ c\ne\res_\t(k)}}
         \frac{L_k-c}{\res_\t(k)-c}$.
\item Let $f_{\s\t}=F_\s m_{\s\t}F_\t$.
\end{enumerate}\end{Definition}

There are some remarks worth making about the definition of $F_\t$.
First, we do not need to specify an order for the product in the
definition of $F_\t$ since the $L_k$ generate a commutative subalgebra
of $\H$. Secondly, the definition of $F_\t$ is very conservative in
the sense that many of the factors of $F_\t$ can be omitted without
changing the element $f_{\s\t}$. Finally, historically this
construction has been used to produce an orthogonal basis for the
Specht modules; we are going to modify this procedure to give an
orthogonal basis for the whole of $\H$.

Let $\s$ an $\t$ be two standard tableaux, not necessarily of the same
shape.  The proof of next result rests upon the easy fact
\cite[Lemma~3.12]{JM:cyc-Schaper} that, {\it because} $P_\H(\Q)$ is
invertible, $\s=\t$ if and only if $\res_\s(k)=\res_\t(k)$ for
$k=1,\dots,n$. 

\begin{Proposition}
Suppose that $P_\H$ is invertible in $R$ and that $\s$ and $\t$
are standard $\lambda$--tableau and
that $k$ is an integer with $1\le k\le n$. Then
\begin{enumerate}
\item 
$f_{\s\t}=m_{\s\t}
             +\Sum_{(\u,\v)\gdom(\s,\t)}a_{\u\v}m_{\u\v}$
for some $a_{\u\v}\in R$;
\item $f_{\s\t}L_k=\res_\t(k)f_{\s\t}$; 
\item $f_{\s\t}F_\u=\delta_{\t\u}f_{\s\t}$; and,
\item $F_\u f_{\s\t}=\delta_{\s\u}f_{\s\t}$.
\end{enumerate}\label{f_st properties}
\end{Proposition}

\begin{proof}
Given (\ref{L_k action}), this is a variation on a well known
argument; see, for example, \cite[Prop. 3.35]{M:ULect}. The proof in
\cite{M:ULect} can be copied out verbatim except that $N$ should be
replaced by $\sum_{k=1}^n|\R(k)|$.
\end{proof}

In particular, part (i) together with (\ref{std basis}) shows that
$$\set{f_{\s\t}|\s, \t\in\Std(\lambda)\ForSome\lambda\vdash n}$$ 
is a basis of~$\H$; shortly we will see that it is an orthogonal basis
of $\H$ with respect to the trace form $\tau$. As a first step we
describe the action of $\H(\Sym_n)$ on this basis; note that the
action of the $L_k$ (and, in particular, $T_0=L_1$), on this basis is
given by Proposition~\ref{f_st properties}(ii).

\begin{Proposition}
Suppose that $\t=\s(i,i+1)$ where $\s$ and $\u$ are standard
$\lambda$--tableaux and $i$ is an integer with $1\le i<n$. If $\t$ is
standard then
\begin{align*}
f_{\u\s}T_i&=\begin{cases}
    \frac{(q-1)\res_\t(i)}{\res_\t(i)-\res_{\s}(i)}f_{\u\s}
       +f_{\u\t},& \If \s\gdom\t,\\[5pt]
    \frac{(q-1)\res_\t(i)}{\res_\t(i)-\res_{\s}(i)}f_{\u\s}
    +\frac{(q\res_\s(i)-\res_\t(i))(\res_\s(i)-q\res_\t(i))}
         {(\res_\t(i)-\res_{\s}(i))^2}f_{\u\t},& \If \t\gdom\s.
   \end{cases}\\
\intertext{If $\t$ is not standard then}
f_{\u\s}T_i&=\begin{cases}
  qf_{\u\s}, &\text{if $i$ and $i+1$ are in the same row of $\s$},\\
  -f_{\u\s}, &\text{if $i$ and $i+1$ are in the same column of $\s$.}
\end{cases}\end{align*}
\label{Tf multiplication}
\end{Proposition}

\begin{proof}
By the above remarks, $\{f_{\u\s}\}$ is a basis of $\H$, so
$f_{\u\s}T_i=\sum_{\a,\b}r_{\a\b}f_{\a\b}$ for some $r_{\a\b}\in R$.
By Proposition~\ref{f_st properties}(iv), $F_\u f_{\a\b}=\delta_{\a\u}f_{\u\b}$.
Therefore, multiplying the equation for $f_{\u\s}T_i$ on the left by
$F_\u$ shows that $r_{\a\b}=0$ whenever $\a\ne\u$; in particular,
$r_{\a\b}=0$ if $\Shape(\b)\ne\lambda$. Hence,  
$f_{\u\s}T_i=\sum_\b a_\b f_{\u\b}$, for some $a_\b\in R$, where $\b$
runs over the set of standard $\lambda$--tableaux. The argument of
\cite[Theorem~3.36]{M:ULect} can now be repeated, essentially word for
word, to complete the proof.
\end{proof}

For each standard $\lambda$--tableau $\s$ let 
$f_\s=f_{\tlam\s}+\Hlam$.  Then $\set{f_\s|\s\in\Std(\lambda)}$
is a basis of the Specht module $S^\lambda$ by 
Proposition~\ref{f_st properties}(i). Note that $f_\tlam=m_\tlam$; further, by
Proposition~\ref{Tf multiplication}, if $\s$ and $\t$ are standard
$\lambda$--tableaux with $\s=\t(i,i+1)\gdom\t$ then
$f_\t=f_\s(T_i-\alpha)$, where
$\alpha=(q-1)\res_\t(i)/\(\res_\t(i)-\res_\s(i)\)$; hence, by
\cite[Cor.~3.14]{JM:cyc-Schaper}, $\set{f_\s|\s\in\Std(\lambda)}$ is
the orthogonal basis of $S^\lambda$ which was constructed in
\cite[Theorem~3.13(vi)]{JM:cyc-Schaper}. (In the notation of
\cite{JM:cyc-Schaper}, $f_\s=\psi_\s(m_\lambda)$.)

The inner products $\<f_\s,f_\t\>$, for $\s,\t\in\Std(\lambda)$ were
computed explicitly (as rational functions) in \cite{JM:cyc-Schaper};
to describe this we need some more notation.  Recall that a node
$y\notin[\lambda]$ is an {\sf addable node} for $\lambda$ if
$[\lambda]\cup\{y\}$ is the diagram of a multipartition; similarly,
$y\in[\lambda]$ is a {\sf removable node} for $\lambda$ if
$[\lambda]\setminus\{y\}$ is the diagram of a multipartition.  Given
two nodes $x=(i,j,k)$ and $y=(a,b,c)$, write $y<x$ if either $c<k$, or
$c=k$ and $b>j$.

Let $\s$ be a $\lambda$--tableau. Then for each integer $i$ there is a
unique node $x\in[\lambda]$ such that $\s(x)=i$, $1\le i\le n$.  Let
$\mathscr A_\s(i)$ be the set of addable nodes for $\Shape(\s\rest i)$
which are strictly less than $x$ (with respect to $<$); similarly, let
$\mathscr R_\s(i)$ be the set of removable nodes strictly less than
$x$ for the multipartition $\Shape(\s\rest i-1)$. Essentially as in
\cite[Defn.~3.15]{JM:cyc-Schaper} let
$$\gamma_\s=q^{\len(d(\s))+\alpha(\lambda)}\prod_{i=1}^n
   \frac{\prod_{y\in\mathscr A_\s(i)}\(\res_\s(i)-\res(y)\)}%
        {\prod_{y\in\mathscr R_\s(i)}\(\res_\s(i)-\res(y)\)}
\Number{gamma definition}$$
where
$\alpha(\lambda)=\frac12\Sum_{s=1}^r
        \sum_{i\ge1}(\lambda^{(s)}_i-1)\lambda^{(s)}_i$. 

If $k$ is an integer let $[k]_q=\frac{q^k-1}{q-1}$ if $q\ne1$; more
generally, set $[k]_q=1+q+\dots+q^{k-1}$ if $k\ge0$ and
$[k]_q=-q^k[-k]_q$ if $k<0$. When $k\ge0$ we also set
$[k]_q^!=[1]_q[2]_q\dots[k]_q$. Finally, if $\lambda$ is a
multipartition let
$[\lambda]_q^!=\prod_{s=1}^r\prod_{i\ge1}[\lambda^{(s)}_i]_q^!$.

\begin{Point}
{\it}{(James--Mathas~\cite[(3.17)--(3.19)]{JM:cyc-Schaper})}
Suppose that $P_\H(\Q)$ is invertible in $R$ and let 
$\t\in\Std(\lambda)$.
Then $\gamma_\t$ is uniquely determined by the two conditions
\begin{enumerate}
\item $\gamma_{\tlam}=
            [\lambda]_q^!\Prod_{1\le s<t\le r}
             \prod_{(i,j)\in[\lambda^{(s)}]}(q^{j-i}Q_s-Q_t)$; and,
\item if $\s=\t(i,i+1)\gdom\t$ then
$\gamma_\t=\frac{(q\res_\s(i)-\res_\t(i))(\res_\s(i)-q\res_\t(i))}
                {(\res_\s(i)-\res_{\t}(i))^2}\gamma_\s$.
\end{enumerate}\noindent
Furthermore, $\<f_\s,f_\t\>=\delta_{\s\t}\gamma_\t$ for all
$\s,\t\in\Std(\lambda)$.
\label{gamma properties}\end{Point}

\begin{Remarks}
\begin{enumerate}
\item In \cite{JM:cyc-Schaper} the rational functions $\gamma_\t$
were important only up to a power of $q$ (which is a unit in $R$).  In
this paper we need to know the inner products $\<f_\s,f_\t\>$ exactly;
for this reason our definition of $\gamma_\s$ differs from that
in~\cite{JM:cyc-Schaper} by the factor
$q^{\len(d(\s))+\alpha(\lambda)}$. The argument in
\cite{JM:cyc-Schaper} computes $\<f_\s,f_\t\>$ using only (\ref{gamma properties})(i) and (\ref{gamma properties})(ii); the inner product
$\<f_\tlam,f_\tlam\>$ is easily seen to be given by the formula on the
right hand side of (\ref{gamma properties})(i). 
\item A factor of $q$ was omitted from the formula for
$\gamma_\t$ in \cite[Lemma~3.17)]{JM:cyc-Schaper}; this is corrected
in (\ref{gamma properties})(ii) --- note that $\len(d(\t))=\len(d(\s))+1$. 
\item The definition of $\gamma_\s$ above is simpler than that given
in \cite{JM:cyc-Schaper} because $\s$ is a standard tableau (rather
than the more general semistandard tableau which were considered
in~\cite{JM:cyc-Schaper} ).
\end{enumerate}\end{Remarks}

Let $(\ ,\ )$ be the inner product on $\H$ given by
$(h_1,h_2)=\tau(h_1h_2^*)$, for $h_1,h_2\in\H$. Then $(\ ,\ )$ is a
symmetric associative bilinear form on $\H$.

\begin{Theorem}
Suppose that $P_\H(\mathbf q)$ is invertible in $R$.
Then
$$\set{f_{\s\t}|\s,\t\in\Std(\lambda)\ForSome\lambda\in\Lambda^+}$$
is a orthogonal basis of $\H$ with respect to the trace
form~$\tau$.  In addition,  if $\s,\t,\u$ and $\v$ are standard
tableaux then $f_{\s\t}f_{\u\v}=\delta_{\u\t}\gamma_\t f_{\s\v}$.  
\label{basis theorem}
\end{Theorem}

\begin{proof}
By (\ref{std basis}) the set $\{m_{\s\t}\}$ is a basis of $\H$;
therefore, $\{f_{\s\t}\}$ is also a basis of~$\H$ by 
Proposition~\ref{f_st properties}(i). 
Next we prove that $f_{\s\t}f_{\u\v}=\delta_{\u\t}\gamma_\t f_{\s\v}$.
First, if $\u\ne\t$ then 
$f_{\s\t}f_{\u\v}=f_{\s\t}F_\u m_{\u\v}F_\v=0$ 
by Proposition~\ref{f_st properties}(iii). Now consider $f_{\s\t}f_{\t\v}$; 
since $\{f_{\a\b}\}$ is a basis we can write
$f_{\s\t}f_{\t\v}=\sum_{\a,\b}r_{\a\b}^\t f_{\a\b}$
for some $r_{\a\b}^\t\in R$. Applying Proposition~\ref{f_st properties}(iii)
and its left handed analogue shows that
$$f_{\s\t}f_{\t\v}=F_\s f_{\s\t}f_{\t\v}F_{\v}
      =\sum_{\a,\b}r_{\a\b}^\t F_\s f_{\a\b}F_\v
      =r_{\s\v}^\t f_{\s\v}.$$
Now, for any $\s,\v\in\Std(\lambda)$ the definition of the inner
product on the Specht module gives
$$\<f_\t,f_\t\>f_{\s\v}\equiv f_{\s\t}f_{\t\v}\mod\Hlam.$$ 
Therefore,
$r_{\s\v}^\t=\<f_\t,f_\t\>=\gamma_\t$ by (\ref{gamma properties}); so
$r_{\s\v}^\t$ depends only on $\t$ and $f_{\s\t}f_{\t\v}=\gamma_\t
f_{\s\v}$ as claimed.

Finally, it remains to show that the basis $\{f_{\s\t}\}$ is
orthogonal with respect to the bilinear form $(\ ,\ )$. First,
$(f_{\s\t},f_{\u\v})=\tau(f_{\s\t}f_{\u\v}^*)=\tau(f_{\s\t}f_{\v\u})
               =\delta_{\t\v}\gamma_\t\tau(f_{\s\u})$;
in particular, $(f_{\s\t},f_{\u\v})=0$ if $\t\ne\v$. On the other
hand, $\tau$ is a trace form so
$(f_{\s\t},f_{\u\v})=\tau(f_{\s\t}f_{\v\u})
                    =\tau(f_{\v\u}f_{\s\t})
                    =\delta_{\s\u}\gamma_\s\tau(f_{\v\t})$.
Therefore, 
$(f_{\s\t},f_{\u\v})=\delta_{\s\u}\delta_{\t\v}\gamma_\t\tau(f_{\s\s})
                    =\delta_{\s\u}\delta_{\t\v}\gamma_\s\tau(f_{\t\t})$.
Consequently, $\{f_{\s\t}\}$ is an orthogonal basis of $\H$ (and
$\tau(f_{\t\t})$ is non--zero for all~$\t$).
\end{proof}

The basis $\{f_{\s\t}\}$ is cellular (with respect to the involution
$*$); but this is not surprising as is was constructed from a cellular
basis.

\begin{Remark}
In \cite{MM:trace} it was shown that if $\H$ is defined over
a ring $R$ in which the parameters $q,Q_1,\dots,Q_r$ are invertible
then $\H$ is a symmetric algebra with respect to the trace
form~$\tau$; however, this was proved indirectly without constructing
a pair of dual bases. The Theorem gives a self--dual basis of the
semisimple Ariki--Koike algebras; no such basis is known in general.
\end{Remark}

As a first consequence, Theorem~\ref{basis theorem} identifies a
submodule of $\H$ which is isomorphic to the Specht module~$S^\lambda$.

\begin{Corollary}
Suppose that $P_\H(\Q)$ is invertible in $R$ and let $\s$ and
$\t$ be standard $\lambda$--tableaux. Then $S^\lambda\cong
f_{\s\t}\H=\sum_{\v\in\Std(\lambda)}Rf_{\s\v}$.
\label{S=fH}
\end{Corollary}

\begin{proof}
By the Theorem, $f_{\s\t}\H$ has as basis the set
$\set{f_{\s\v}|\v\in\Std(\lambda)}$. The isomorphism is
given by the linear map $S^\lambda\To f_{\s\t}\H$ determined by
$f_\v\mapsto f_{\s\v}$ for all tableaux $\v\in\Std(\lambda)$.
\end{proof}

Set $\tilde f_{\s\t}=\gamma_\t^{-1}f_{\s\t}$.  Then 
$\tilde f_{\s\t}\tilde f_{\u\v}=\delta_{\t\u}\tilde f_{\s\v}$ and
$\{\tilde f_{\s\t}\}$ is a basis of $\H$. Hence, we have the
following. 

\begin{Corollary}
Suppose that $P_\H(\Q)$ is invertible in $R$. Then
$$\set{\tilde f_{\s\t}|\s,\t\in\Std(\lambda)\ForSome\lambda\in\Lambda^+}$$
is a basis of matrix units in $\H$.
\end{Corollary}

The last result yields an explicit isomorphism from $\H$ to the group
ring of $W_{r,n}$ when~$P_\H(\Q)$ is invertible.  Assume that $R$
contains a primitive $r$th root of unity $\zeta$; then, $\H\cong
RW_{r,n}$ when $q=1$ and $Q_s=\zeta^s$ for $s=1,2,\dots,r$. Write
$\tilde f^1_{\s\t}$ for the element of~$RW_{r,n}$ corresponding to
$\tilde f_{\s\t}\in\H$ under the canonical isomorphism
$\H\To RW_{r,n}$.

\begin{Corollary}
Assume that $R$ contains a primitive $r$th root of unity and that
$P_\H(\Q)$ is invertible in $R$. Then $\H\cong R W_{r,n}$
via the $R$--algebra homomorphism determined 
by~$\tilde f_{\s\t}\mapsto \tilde f^1_{\s\t}$.
\label{isomorphism}\end{Corollary}

By parts (i) and (iii) of Theorem~\ref{idempotents} below,
$T_i=\sum_\t \tilde f_{\t\t} T_i$, for $0\le i<n$;
so, in principle, we can determine the image of the generators of
$\H$ under this isomorphism.

Lusztig~\cite{L:iso} has shown that there exists a homomorphism $\Phi$
from the Hecke algebra $\H(W)$ of any finite Weyl group $W$ to the
group ring $RW$ and he shows that $\Phi$ induces an isomorphism when
$\H(W)$ is semisimple. Our map is not an analogue of Lusztig's
isomorphism; rather it is an explicit realization of the Tits
deformation theorem in this setting. It would be good to find a
generalization of Lusztig's isomorphism theorem for the Ariki--Koike
algebras.

We next construct the primitive (central) idempotents in $\H$. 

\begin{Theorem}
Suppose that $R$ is a field and that $P_\H(\Q)\ne0$.
\begin{enumerate}
\item Let $\t$ be a standard $\lambda$--tableau. Then
$F_\t=\frac1{\gamma_\t}f_{\t\t}$ and $F_\t$ is a primitive idempotent
with $S^\lambda\cong F_\t\H$.
\item  For any multipartition $\lambda$ let 
$F_\lambda=\sum_{\t\in\Std(\lambda)} F_\t$. Then
$F_\lambda$ is a primitive central idempotent.
\item $\set{F_\lambda|\lambda\vdash n}$ is a complete set of
primitive central idempotents; in particular,
$$1=\sum_{\lambda\vdash n}F_\lambda
   =\sum_{\t\standard}F_\t.$$
\end{enumerate}\label{idempotents}
\end{Theorem}

\begin{proof}
We may write $F_\t=\sum_{\u,\v}a_{\u\v}f_{\u\v}$ for some
$a_{\u\v}\in R$ by Theorem~\ref{basis theorem}. By 
Proposition~\ref{f_st properties}(iii) and Theorem~\ref{basis theorem},
$f_{\s\t}=f_{\s\t}F_\t=\sum_{\u,\v}a_{\u\v}f_{\s\t}f_{\u\v}
          =\sum_\v \gamma_\t a_{\t\v} f_{\s\v};$
equating coefficients on both sides shows that $a_{\t\v}=0$ if
$\v\ne\t$ and that $1=a_{\t\t}\gamma_t$. Since $F_\t^*=F_t\PS$ we also
have that $a_{\v\t}=0$ if $\v\ne\t$. Hence,
$F_\t=\frac1{\gamma_\t}f_{\t\t}$ as claimed.  By Theorem~\ref{basis theorem},
$\frac1{\gamma_\t}f_{\t\t}$ is idempotent; further, it is primitive
because~$S^\lambda$ is irreducible and $S^\lambda\cong
\frac1{\gamma_\t}f_{\t\t}\H$ by Corollary~\ref{S=fH}. Hence, (i) is proved.

Parts (ii) and (iii) now follow because 
$\H=\bigoplus_{\lambda\vdash n}\bigoplus_{\t\in\Std(\lambda)}F_\t\H$ 
is a decomposition of $\H$ into a direct sum of simple modules with
each simple module $F_\t\H\cong S^\lambda$ appearing with
multiplicity equal to its dimension (the sum is direct
because $\{f_{\t\v}\}$ is a basis of $F_\t\H$ by Theorem~\ref{basis theorem}
and the set of all $f_{\u\v}$ is a basis of~$\H$).
\end{proof}

\begin{Corollary}
Let $\t$ be a standard tableau and let $k$ be an
integer with $1\le k\le n$. Then $F_\t L_k=\res_\t(k)F_\t$
and $f_{\s\t}L_k=\res_\t(k)f_{\s\t}$.
\label{F_t L_k} 
\end{Corollary}

\begin{proof}
Since $F_\t=\frac 1{\gamma_\t}f_{\t\t}$ by part (i), the
formula for $F_\t L_k$ follows from Proposition~\ref{f_st properties}(ii); this
also implies the second statement because $f_{\s\t}L_k=f_{\s\t}F_\t L_k$.
\end{proof}

In particular, the Corollary describes the action of $T_0=L_1$ on
the orthogonal basis $\{f_{\s\t}\}$ of $\H$. Proposition~\ref{Tf multiplication}
gives the action of the remaining generators~$T_1,\dots,T_{n-1}$ 
of~$\H$ on the basis $\{f_{\s\t}\}$.

\begin{Corollary}
Suppose that $1\le k\le n$. Then $\prod_{c\in\R(k)}(L_k-c)=0$ and
this is the minimum polynomial for $L_k$ acting on $\H$.
\label{min poly}
\end{Corollary}

\begin{proof}
By Theorem~\ref{idempotents},
$$\prod_{c\in\R(k)}(L_k-c)
  =1\cdot\prod_{c\in\R(k)}(L_k-c)
  =\sum_{\t\text{\ standard}}F_\t\prod_{c\in\R(k)}(L_k-c)=0,$$
since $F_\t\(L_k-\res_\t(k)\)=0$ by Corollary~\ref{F_t L_k}.
Moreover, if we remove any factor $(L_k-c)$ from the product
$\prod_{c\in\R(k)}(L_k-c)$ then what remains is non--zero because it
divides $F_\t$ for some standard tableau $\t$. Hence,
$\prod_{c\in\R(k)}(L_k-c)$ is the minimum polynomial of $L_k$.
\end{proof}

\begin{Corollary}
Suppose that $1\le k\le n$. Then
$L_k=\sum_\t \res_\t(k) F_\t$, where the sum is over the set of all
standard tableaux $($of arbitrary shape$)$.
\label{L-F_t lemma}
  
\end{Corollary}

\begin{proof}
Combining part (iii) of Theorem~\ref{idempotents} with Corollary~\ref{F_t L_k}
shows that
$$L_k=1\cdot L_k
     =\sum_{\t\standard} F_\t L_k
     =\sum_{\t\standard} \res_\t(k)F_\t$$
as required.
\end{proof}

We close this section with a description of the centre of $\H$ in the
semisimple case; this result is due to
Ariki and Koike~\cite[Theorem~3.20]{AK1}.

\begin{Theorem}
[Ariki--Koike] Suppose that $R$ is a field and that 
$P_\H(\Q)\ne0$.
\begin{enumerate}
\item The centre of $\H$ is the set of symmetric
polynomials in $L_1,L_2,\dots,L_n$.
\item Let $\mathbb L$ be the subalgebra of $\H$ generated by 
$L_1,L_2,\dots,L_n$. Then $\mathbb L$ is a maximal
abelian subalgebra of $\H$.
\end{enumerate}
\end{Theorem}

\begin{proof}
First consider (ii). By definition $\mathbb L$ contains each of
the primitive idempotents~$F_\t$, for an arbitrary standard
tableau~$\t$. On the other hand, by Corollary~\ref{L-F_t lemma},~$\mathbb L$ is
contained in the subalgebra of $\H$ generated by the primitive
idempotents $F_\t$; hence,~$\mathbb L$ is the subalgebra of $\H$
generated by the idempotents~$F_\t$. As the primitive idempotents
generate a maximal abelian subalgebra of $\H$ the result follows.

Now consider~(i). By (\ref{L-comm}) every symmetric polynomial in
$L_1,\dots,L_n$ belongs to the centre of~$\H$. Conversely, the centre
of $\H$ has as basis the set of idempotents $F_\lambda$, as $\lambda$
runs over the multipartitions of $n$. So to prove (i) it is enough to
show that each $F_\lambda$ is symmetric in $L_1,\dots,L_n$. 

Let $\R=\cup_{k=1}^n\R(k)$ be the set of all possible residues for
$\H$. By Corollary~\ref{F_t L_k} if  $\t$ is a $\lambda$--tableau and 
$1\le m\le n$ then $F_\t(L_k-c)=(\res_\t(k)-c)F_\t$; therefore,
multiplying $F_\t$ by the appropriate extra factors we see that
$$ F_\t\defeq\prod_{k=1}^n\prod_{c\in\R(k)\setminus\{\res_\t(k)\}}
            \frac{L_k-c}{\res_{\t}(k)-c}
       =\prod_{k=1}^n\prod_{c\in\R\setminus\{\res_{\t}(k)\}}
            \frac{L_k-c}{\res_{\t}(k)-c}.$$
Notice that the denominator is equal to
$f_\lambda=\prod_{k=1}^n\prod_{c\in\R\setminus\{\res_{\tlam}(k)\}}
  (\res_{\tlam}(k)-c)$
and $f_\lambda$ depends only on~$\lambda$ and not directly on $\t$.

Next, following Murphy~\cite{murphy:hecke} say that a
$\lambda$--tableau is {\sf regular} if its entries increase from left to
right along the nodes in $[\lambda]$ of constant residue. Now, because
$P_\H(\Q)\ne0$ two nodes in $[\lambda]$ have the same residue if and
only if they lie on the same diagonal
$\set{(i+d,j+d,s)\in[\lambda^{(s)}]|d\ge0}$ of $[\lambda]$; thus, a
tableau is regular if and only if its entries increase from left to
right along each diagonal. (So, for example, every standard tableau is
regular, but not conversely.) Extending the formula above for
$F_\t$, for each regular tableau~$\t$ define
$$ F_\t=\prod_{k=1}^n\prod_{c\in\R\setminus\{\res_\t(k)\}}
            \frac{L_k-c}{\res_{\t}(k)-c}
       =\frac1{f_\lambda}
  \prod_{k=1}^n\prod_{c\in\R\setminus\{\res_\t(k)\}}(L_k-c).$$
Observe that if $c\in\R$ then $L_k-c$ is a factor of $F_\t$ if and
only if $c\ne\res_\t(k)$; therefore,~$F_\t$ determines $\res_\t(k)$,
for $k=1,\dots,n$. As remarked above, the residues on
the different diagonals of~$[\lambda]$ are distinct; consequently, a
regular tableau $\t$ is uniquely determined by the sequence of
residues $\(\res_\t(1),\dots,\res_\t(n)\)$ and hence by the
`polynomial'~$F_\t$. It follows that if we permute $L_1,\dots,L_n$
then $F_\t$ is mapped to~$F_\s$ where $\s$ is the regular tableau
determined by the corresponding permutation of the residue sequence
of~$\t$~--- again,  $\s$ is necessarily a $\lambda$--tableau because
the shape of $\s$ is determined by the lengths of its diagonals which,
in turn, are determined by the multiplicity of each residue in $\s$
(or $\t$).

Finally, notice that a regular tableau $\t$ is {\it not} standard
if and only if $\res_\t(k)\notin\R(k)$ for some $k$. Therefore, if
$\t$ is not standard then $\prod_{c\in\R(k)}(L_k-c)$ is a factor of
$F_\t$ and, consequently, $F_\t=0$ by Corollary~\ref{min poly}; hence,
$$F_\lambda
 =\sum_{\substack{\t\text{\ a standard}\\\lambda\text{--tableau}}}F_\t
 =\sum_{\substack{\t\text{\ a regular}\\\lambda\text{--tableau}}}F_\t.$$
By the last paragraph the right hand side is a symmetric polynomial in
$L_1,\dots,L_n$ so the theorem follows.
\end{proof}

For the Iwahori--Hecke algebras of type $A$ (that is, when $r=1$),
it is conjectured that the centre of $\H$ is always the set of
symmetric polynomials in $L_1,\dots,L_n$, even when $\H$ is not
semisimple; see \cite{DJ:blocks,M:centre}. When $r>1$ and $\H$ is
not semisimple there are cases where the centre of $\H$ is
larger than the set of symmetric polynomials in~$L_1,\dots,L_n$;
for an example see \cite[p.~792]{Ariki:can}.

\section{Another construction of the Specht modules}

In this section we give another two constructions of the Specht
modules.  The first is via a second cellular basis of $\H$ which is,
in a certain sense, dual to the basis described in the previous
section. The second construction combines these two approaches to
produce submodules of $\H$ which are isomorphic to the Specht modules
$S^\lambda$. Some of these results can be found in the work of Du and
Rui~\cite{DuRui:branching}. In the next section we will use these
results to compute the Schur elements of~$\H$.

Let $\ZZ=\Z[\hat q,\hat q^{-1},\hat Q_1,\dots,\hat Q_r]$, where
$\hat q, \hat Q_1,\dots,\hat Q_r$ are indeterminates over $\Z$,
and let $\HZ$ be the Ariki--Koike algebra with parameters 
$\hat q,\hat Q_1,\dots,\hat Q_r$. Consider the ring $R$ as a
$\ZZ$--module by letting $\hat q$ act on $R$ as multiplication by
$q$ and $\hat Q_s$ by multiplication by $Q_s$, for 
$1\le s\le r$. Then $\H\cong\HZ\otimes_\ZZ R$, since $\H$ is free as
an $R$--module; we say that $\H$ is a {\sf specialization} of $\HZ$ and
call the map which sends $h\in\HZ$ to $h\otimes1\in\H$ the
{\sf specialization homomorphism}.  

Let ${}'\map\ZZ\ZZ$ be the $\Z$--linear map given by 
$\hat q\mapsto\hat q^{-1}$ and $\hat Q_s\mapsto \hat Q_{r-s+1}$ for
$1\le s\le r$. Define $T_0'=T_0$ and $T_i'=-\hat q^{-1}T_i$ for 
$1\le i<n$; using the relations of $\HZ$ it is easy to verify
that $'$ now extends to a $\Z$--linear ring involution $'\map\HZ\HZ$
of $\HZ$. Hereafter, we drop the distinction between $\hat q$ and $q$,
and $\hat Q_s$ and~$Q_s$.

Suppose that $h\in\H$. Then there exists a (not necessarily unique)
$h_\ZZ\in\HZ$ such that $h=h_\ZZ\otimes 1$ under specialization; we
sometimes abuse notation and write $h'=h_\ZZ'\otimes 1\in\H$. As the
map ${}'$ does not in general define a semilinear involution on $R$,
this notation is not well--defined on elements of $\H$; however, in
the cases where we employ it there should be no ambiguity. For
example, $T_w'=(-q)^{-\len(w)}T_w$ and $L_i'=L_i$ for all $w\in\Sym_n$
and $1\le i\le n$.

If $\lambda$ is a multipartition of $n$ let
$y_\lambda=\sum_{w\in\Sym_\lambda}(-q)^{-\len(w)}T_w$; then
$y_\lambda=x_\lambda'$. Similarly, we define 
$n_\lambda=y_\lambda u_\lambda^-=m_\lambda'$ where
$$u_\lambda^-=\prod_{s=2}^r\prod_{k=1}^{\a_s}(L_k-Q_{r-s+1})
             =\prod_{s=1}^{r-1}\prod_{k=1}^{\a_{r-s+1}}(L_k-Q_{s}).$$
Observe that $u_\lambda^-=(u_\lambda^+)'$; here, as usual,
$\a_s=|\lambda^{(1)}|+\dots+|\lambda^{(s-1)}|$ for all $s$. For
standard tableaux $\s,\t\in\Std(\lambda)$ set
$n_{\s\t}={T_{d(\s)}^*}'n_\lambda{T_{d(\t)}^{\phantom*}}'$; then
$m_{\s\t}'\in\HZ$ is mapped to $n_{\s\t}$ under specialization.
Hence, from (\ref{std basis}) we obtain the following.

\begin{Point}
{\it}{\cite[(2.7)]{DuRui:branching}}
The Ariki--Koike algebra $\H$ is free as an $R$--module with
cellular basis
$\set{n_{\s\t}|\s,\t\in\Std(\lambda)\ForSome\lambda\in\Lambda^+}$.
\end{Point}

Let $\lambda$ be a multipartition of $n$. Then $(\Hlam)'$ is a
two--sided ideal of $\H$ which is free as an $R$--module with basis
$\set{n_{\u\v}|\u,\v\in\Std(\mu)\ForSome\mu\gdom\lambda}$. 
Let~$\tilde S^\lambda$ be the Specht module (or cell module) 
corresponding to $\lambda$ determined by the basis $\{n_{\s\t}\}$; 
then $\tilde S^\lambda\cong n_\lambda\H/(n_\lambda\H\cap(\Hlam)')$
and $\tilde S^\lambda$ is free as an $R$--module with basis
$\set{n_\t|\t\in\Std(\lambda)}$, where
$n_\t=n_{\tlam\t}+(\Hlam)'=\(m_{\tlam\t}+\Hlam)'$ for all
$\t\in\Std(\lambda)$.

In order to compare the two modules $S^\lambda$ and $\tilde
S^\lambda$ we need to introduce some more notation.
Given a partition $\sigma$ let $\sigma'=(\sigma_1',\sigma_2',\dots)$
be the partition which is conjugate to~$\sigma$; thus, $\sigma_i'$
is the number of nodes in column $i$ of the diagram of $\sigma$.
If $\lambda=\rtuple\lambda$ is a multipartition then the {\sf conjugate}
$\lamp=\rtuple{\lamp}$ of $\lamp$ is the multipartition
with $\lamp^{(s)}=(\lambda^{(r-s+1)})'$ for $1\le s\le r$.

Now suppose that $\t=\rtuple\t$ is a standard $\lambda$--tableau.
Then the {\sf conjugate} of~$\t$ is the standard $\lamp$--tableau
$\t'=\rtuple{{\t'}}$ where ${\t'}^{(s)}$ is the tableau obtained by
interchanging the rows and columns of $\t^{(r-s+1)}$.

With these definitions in place, we see that the following holds in
$\ZZ$.

\begin{Point}
{\it}* Let $\t$ be a standard $\lambda$--tableau. Then 
$\(\res_\t(k)\)'=\res_{\t'}(k)$ in $\ZZ$, for $1\le k\le n$.
\label{res'}\end{Point}

The expression $\res_{\t'}(k)$ is always well--defined; whereas
$(\res_\t(k))'$ is ambiguous for certain rings $R$. As a first
consequence we determine how the $L_k$ act on the $n_{\s\t}$--basis 
of~$\H$.

\begin{Proposition}
Let $\s$ and $\t$ be standard $\lambda$--tableaux and suppose
that $k$ is an integer with $1\le m\le n$. Then there exist 
$a_\v\in R$ such that
$$n_{\s\t}L_k=\res_{\t'}(k)n_{\s\t}
       +\sum_{\SR{\v\in\Std(\lambda)}{\v\gdom\t}}a_\v n_{\s\v}
                \bmod(\Hlam)'.$$
\label{L_k' action}
\end{Proposition}

\begin{proof}
First assume that $R=\ZZ$. Then ${}'$ is a $\Z$--linear ring
involution on $\HZ$ and~$L_k'=L_k$; therefore, by~(\ref{L_k action}),
$$n_{\s\t}L_k=(m_{\s\t}L_k)'
     =\Big(\res_{\t}(k)m_{\s\t} +\sum_{\SR{\v\in\Std(\lambda)}{\v\gdom\t}}a_\v
       m_{\s\v}\bmod\Hlam\Big)'.$$
Using (\ref{res'}) this proves the Proposition for $\HZ$. The general
case now follows by specialization since $\H\cong\HZ\otimes_\ZZ R$.
\end{proof}

Next consider the orthogonal basis $\{f_{\s\t}\}$ of $\H$ in the
case where $P_\H(\Q)$ is invertible. Let $\ZZ_P$ be the localization
of $\ZZ$ at~$P_\H(\Q)$ and let~$\H_{\ZZ_P}$ be the corresponding
Ariki--Koike algebra.  The involution ${}'$ extends to~$\H_{\ZZ_P}$
and~$\H$ is a specialization of~$\H_{\ZZ_P}$ whenever $P_\H(\Q)$ is
invertible in $R$. (Note that $Q_1,\dots,Q_r$ are indeterminates in
$\ZZ_P$.) 

In general, $f_{\s\t}\notin\HZ$; however, if
$\t\ne\u$ then $\res_\t(k)-\res_\u(k)$ is a factor of~$P_\H(\Q)$ for
all $k$, so $f_{\s\t}\in\H_{\ZZ_P}$ and we can speak of the elements
$F_\t$ and $f_{\s\t}\in\H_{\ZZ_P}$. More generally, whenever
$P_\H(\Q)$ is invertible in $R$ we have an element $f_{\s\t}'\in\H$
via specialization because $\H\cong\H_{\ZZ_P}\otimes_{\ZZ_P}R$.

\begin{Proposition}
Suppose that $\t$ is a standard tableau. Then $F_\t'=F_{\t'}$
in $\H_{\ZZ_P}$.
\label{F'}
\end{Proposition}

\begin{proof}
Applying the definitions together with (\ref{res'}) gives
\begin{align*}
F_\t' &=\prod_{k=1}^n
   \prod_{\substack{c\in\R(k)\\c\ne\res_\t(k)}}
 \Big(\frac{L_k-c}{\res_\t(k)-c}\Big)'
       =\prod_{k=1}^n\prod_{\substack{c\in\R(k)\\c\ne\res_{\t}(k)}}
  \frac{L_k-c'}{\res_{\t'}(k)-c'}\\
      &=\prod_{k=1}^n\prod_{\substack{c\in\R(k)\\c\ne\res_{\t'}(k)}}
  \frac{L_k-c}{\res_{\t'}(k)-c}
       =F_{\t'},
\end{align*}
the last equality following because $\R(k)$ is invariant under ${}'$.
\end{proof}

By Lemma~\ref{generic degree lemma} and Theorem~\ref{idempotents} the Schur elements
are given by $s_\lambda(\Q)=\tau(F_\tlam)^{-1}$; consequently, the
Schur elements have the following ``palindromy'' property.

\begin{Corollary}
Suppose that $\lambda$ is a multipartition of $n$. Then
$s_\lamp(\Q)=\(s_\lambda(\Q)\)'$.
\label{symmetry}
\end{Corollary}

Returning to the general case, let $g_{\s\t}=F_{\s'}n_{\s\t}F_{\t'}$;
then $g_{\s\t}=F_\s'm_{\s\t}'F_\t'=f_{\s\t}'$ (in~$\H_{\ZZ_p}$).
Applying ${}'$ to $\{f_{\s\t}\}$ and using Theorem~\ref{basis theorem} (and a
specialization argument) shows that
$\set{g_{\s\t}|\s,\t\in\Std(\lambda)\ForSome\lambda\in\Lambda^+}$ is a
basis of $\H$. Consequently, as in Corollary~\ref{S=fH}, 
$\tilde S^\lambda\cong
g_{\s\t}\H$ for any standard $\lambda$--tableaux
$\s,\t\in\Std(\lambda)$.

\begin{Remark}
By the Proposition and Theorem~\ref{idempotents}(i),
$$g_{\t\t}=f_{\t\t}'=(\gamma_\t F_\t)'=\gamma_\t'F_{\t'}
         =\frac {\gamma_\t'}{\gamma_{\t'}}f_{\t'\t'}.$$
More generally, we can write $g_{\s\t}=\sum_{\u,\v}a_{\u\v}f_{\u\v}$
for some $a_{\u\v}\in R$. By Proposition~\ref{f_st properties} and Proposition~\ref{F'},
$F_{\s'}g_{\s\t}F_{\t'}=g_{\s\t}$; so it follows that $a_{\u\v}=0$
unless $\u=\s'$ and $\v=\t'$. Therefore,
$g_{\s\t}=\alpha_{\s\t}f_{\s'\t'}$ for some $\alpha_{\s\t}\in R$.
Applying the $*$--involution shows that $\alpha_{\s\t}=\alpha_{\t\s}$.
Finally, by looking at the product $g_{\s\t}g_{\t\s}$ we see that
$\alpha_{\s\t}^2=\gamma_\s'\gamma_\t'/\gamma_{\s'}\gamma_{\t'}$.  
{\it A~priori}, there is no reason why the square root of this element
should belong to~$R$; nor do I see a way to determine the sign of
$\alpha_{\s\t}$.
\label{g_st to f_st}\end{Remark}

Combining Proposition~\ref{F'} with Corollary~\ref{S=fH} and the corresponding result for the
$g$--basis shows that
$S^\lambda\cong f_{\t\t}\H=g_{\t'\t'}\H\cong\tilde S^\lamp$, for any
$\t\in\Std(\lambda)$. Hence, we have the following.

\begin{Corollary}
Suppose that $P_\H(\Q)$ is invertible in $R$. Then 
$\tilde S^\lambda\cong S^\lamp$.
\label{dual isomorphism}\end{Corollary}

When $R$ is field the assumption that $P_\H(\Q)$ is invertible is
equivalent to $\H$ being semisimple. This assumption is necessary
because, in general, $S^\lamp$ and $\tilde S^\lambda$ are not
isomorphic; rather, $S^\lamp$ is isomorphic to the {\it dual} of
$\tilde S^\lambda$~\cite{M:tilting}. In the semisimple case both
$S^\lambda$ and $\tilde S^\lamp$ are irreducible, and hence
self--dual, since they carry a non--degenerate bilinear form.
Accordingly, we call the module~$\tilde S^\lambda$ a {\sf dual Specht
module}.

Here is another useful application of Proposition~\ref{F'}.

\begin{Corollary}
Suppose that $\s$ and $\t$ are standard $\lambda$--tableaux and that
$\u$ and $\v$ are standard $\mu$--tableaux where $\lambda$ and $\mu$ are
multipartitions of $n$. Then $f_{\s\t}g_{\u\v}=0$ if $\t\ne\u'$.
\label{cancellation}
\end{Corollary}

\begin{proof}
Applying the definitions, 
$f_{\s\t}g_{\u\v}
   =F_\s m_{\s\t}F_\t F_\u'm_{\u\v}F_\v'$;
however, $F_\u'=F_{\u'}$ by Proposition~\ref{F'}, so $F_\t
F_\u'=\delta_{\t\u'}F_\t$ by Theorem~\ref{idempotents}(i), giving the
result. (By Remark~\ref{g_st to f_st}, $f_{\s\t}g_{\t'\v}$ is a 
scalar multiple of $f_{\s\v'}$.)
\end{proof}

The Specht modules $S^\lambda$ and the dual Specht modules $\tilde
S^\lambda$ are both constructed as quotient modules using the bases
$\{m_{\s\t}\}$ and $\{n_{\s\t}\}$ respectively (see Corollary~\ref{S=fH}). In
the cases where $P_\H(\Q)$ is invertible in $R$ we have also
constructed these modules as submodules of $\H$. Next we produce
submodules of an arbitrary Hecke algebra~$\H$ which are isomorphic
to the Specht modules; these results will also play a r\^ole in
computing the Schur elements in the next section.

Recall that $\tlam$ is the $\lambda$--tableau which has the numbers
$1,2,\dots,n$ entered in order first along the rows of $\tlam^{(1)}$
and then the rows of $\tlam^{(2)}$ and so on. Let
$\tllam=(\t^\lamp)'$; that is to say that $\tllam$ is the
$\lambda$--tableau with the numbers $1,2,\dots,n$ entered in order
first down the columns of $\tllam^{(r)}$ and then the columns of
$\tllam^{(r-1)}$ etcetera. Observe that if $\t$ is a standard
$\lambda$--tableau then $\tlam\gedom\t\gedom\tllam$.

\begin{Proposition}
Suppose that $P_\H(\Q)$ is invertible and let $\lambda$ be a
multipartition of~$n$. Then $m_\lambda\H n_\lamp=Rf_{\tlam\tllam}$.
\label{one dimensional}
\end{Proposition}

\begin{proof}
By Proposition~\ref{f_st properties}(i),
$m_\lambda=f_{\tlam\tlam}+\sum_{\u,\v\gdom\tlam}a_{\u\v}f_{\u\v}$ for
some $a_{\u\v}\in R$. Therefore, interchanging the roles of $\lambda$
and $\lamp$ and applying the involution ${}'$ (in~$\H_{\ZZ_p}$ and
then specializing) we see that there exist $b_{\a\b}\in R$ such that
$$n_\lamp=g_{\tlamp\tlamp}
             +\sum_{\a,\b\gdom\tlamp}b_{\a\b}g_{\a\b}
      =\frac{\gamma_{\tlamp}'}{\gamma_{\tllam}}f_{\tllam\tllam}
             +\sum_{\tllam\gdom\a',\b'}b_{\a\b}g_{\a\b}
\Number{n_lam'}$$
where for the second equality we have used Remark~\ref{g_st to f_st} (note
that $(\tlamp)'=\tllam$) and the observation that
$\a,\b\gdom\tlamp$ if and only if~$\tllam\gdom\a',\b'$.  Now
$m_\lambda\H n_\lamp$ is spanned by the elements $m_\lambda f_{\s\t}
n_\lamp$, where $\s$ and $\t$ range over all pairs of standard
tableaux of the same shape. Now, (\ref{n_lam'}) and 
Corollary~\ref{cancellation} imply that
\begin{align*}
m_\lambda f_{\s\t} n_\lamp
   &=\Big(f_{\tlam\tlam}+\sum_{\u,\v\gdom\tlam}a_{\u\v}f_{\u\v}\Big)
       f_{\s\t}
     \Big(\frac{\gamma_{\tlamp}'}{\gamma_{\tllam}}f_{\tllam\tllam}
             +\sum_{\tllam\gdom\a',\b'}b_{\a\b}g_{\a\b}\Big)\\
  &=\frac{\gamma_{\tlamp}'}{\gamma_{\tllam}}
       f_{\tlam\tlam}f_{\s\t}f_{\tllam\tllam}
        +\sum_{\substack{\tllam\gdom\a',\b'\\\u,\v\gdom\tlam}}
           a_{\u\v}b_{\a\b}f_{\u\v}f_{\s\t}g_{\a\b},\\
   &=\frac{\gamma_{\tlamp}'}{\gamma_{\tllam}}
                     f_{\tlam\tlam}f_{\s\t}f_{\tllam\tllam}\\
   &=\begin{cases} \gamma_\tlam\gamma_{\tlamp}'
	             f_{\tlam\tllam}, &\If\s=\tlam\And\t=\tllam,\\
	        0,&\Otherwise,
         \end{cases}
\end{align*}
with the last equality following from Theorem~\ref{basis theorem}.
Therefore, $m_\lambda\H n_\lamp=R f_{\tlam\tllam}$ as required.
\end{proof}

Let $w_\lambda=d(\tllam)$; thus, $w_\lambda$ is the unique element of
$\Sym_n$ such that $\tllam=\tlam w_\lambda$.

\begin{Example}
Let $\lambda=\((2,1^2),(2,1),(2)\)$. Then
$$\tlam =\bigg(\begin{array}[b]{*4c}\young(12,3,4)\ ,
             &\begin{array}{@{}c@{}}\young(56,7)\\\strut\end{array} ,
	     &\begin{array}{@{}c@{}}\young(89)\\[3.5mm]\strut\end{array}\ \bigg)\end{array}
\quad\And\quad
\tllam=\bigg(\begin{array}[b]{*3c}\young(69,7,8)\ ,%
                &\begin{array}{@{}c@{}}\young(35,4)\\\strut\end{array} ,%
	        &\begin{array}{@{}c@{}}\young(12)\\[3.5mm]\strut
             \end{array}\ %
       \end{array}\bigg)
$$
and $w_\lambda=(1,6,5,3,7,4,8)(2,9)$.
\label{w_lambda example}\end{Example}

In order to compute the Schur elements we will need to know quite a
few properties of the permutations $w_\lambda$; these permutations
enter the story through the following definition and Corollary.

\begin{Definition}
Suppose that $\lambda$ is a multipartition of $n$. Let
$z_\lambda=m_\lambda T_{w_\lambda} n_\lamp$.
\end{Definition}

The element $z_\lambda$ and the following result are crucial to our 
computation of the Schur elements.

\begin{Corollary}
Suppose that $P_\H(\Q)$ is invertible in $R$. Then 
$z_\lambda=\gamma_{\tlamp}'f_{\tlam\tllam}$. In particular,
$m_\lambda\H n_\lamp=Rz_\lambda$.
\label{z=f}
\end{Corollary}

\begin{proof}
Consulting the definitions,
$z_\lambda=m_\lambda T_{w_\lambda}n_\lamp
          =m_{\tlam\tllam}n_\lamp$. Furthermore,  
by Proposition~\ref{f_st properties}(i), there exist $c_{\u\v}\in R$ such that
$m_{\tlam\tllam}=f_{\tlam\tllam}+\sum c_{\u\v}f_{\u\v}$
where the sum is over the pairs $(\u,\v)$ of standard tableaux which
strictly dominate~$\tllam$. Therefore, by (\ref{n_lam'}),
$$z_\lambda=\Big(f_{\tlam\tllam}
        +\sum_{\u,\v\gdom\tllam}c_{\u\v}f_{\u\v}\Big)
     \Big(\frac{\gamma_{\tlamp}'}{\gamma_{\tllam}}f_{\tllam\tllam}
    +\sum_{\SR{\a,\b}{\tllam\gdom\a',\b'}}b_{\a\b}g_{\a\b}\Big).
$$
By Theorem~\ref{basis theorem},
$f_{\tlam\tllam}f_{\tllam\tllam}=\gamma_{\tllam}f_{\tlam\tllam}$ and
this is the only non--zero term in this product by Corollary~\ref{cancellation}.
Hence, $z_\lambda=\gamma_{\tlamp}'f_{\tlam\tllam}$ as required.
\end{proof}	

Now, $z_\lambda$ is an element of $\HZ$, so
$\gamma_{\tlamp}'f_{\tlam\tllam}\in\HZ$. By definition, $z_\lambda\H$
is a submodule of $m_\lambda\H$ and a quotient module of $n_\lamp\H$.
Over an arbitrary ring $R$, Du and
Rui~\cite[Remark~2.5]{DuRui:branching} showed that $S^\lambda\cong
z_\lambda^*\H$ and $\tilde S^{\lambda'}\cong z_\lambda\H$ as
$\H$--modules, the isomorphisms being given by the natural quotient
maps $m_\lambda\H\To z_\lambda^*\H$ and $n_{\lamp}\H\To z_\lambda\H$.
Note that $S^\lambda\cong\tilde S^{\lamp}$ when $\H$ is semisimple by
Corollary~\ref{S=fH}.

\section{The Schur elements}

Using the results of the previous sections we are now ready to compute
the Schur elements of $\H$. By Lemma~\ref{generic degree lemma},
$s_\lambda(\Q)=\frac1{\tau(F_\tlam)}$, so it is enough to calculate
$\tau(F_\tlam)$. Our basic strategy, which is inspired by
Murphy~\cite{murphy:hecke}, is to write~$F_\tlam$ as a product of two
terms and, in effect, to evaluate $\tau$ on each of these factors
separately. 

\begin{Proposition}
Suppose that $\t$ is a standard $\lambda$--tableau. Then
there exist elements~$\Phi_\t$ and $\Psi_\t$ in $\H(\Sym_n)$ such
that
\begin{enumerate}
\item $\Psi_\t F_\t=F_\tlam\Phi_\t$;
\item $\Psi_\t=T_{d(\t)}+\Sum_{w<d(\t)} p_{\t w}T_w$, for some 
$p_{\t w}\in R$; and,
\item $\glam\Phi_\t\PS\Psi_\t^*=\gamma_\t$.
\end{enumerate}\label{Phi-Psi prop}
\end{Proposition}

\begin{proof}
We prove all three statements by induction on~$\t$. When
$\t=\tlam$ there is nothing to prove as we may take
$\Phi_\tlam=\Psi_\tlam=1$. Suppose then that $\t\ne\tlam$. Then there
exists an integer $i$, with $1\le i<n$, such that
$\s=\t(i,i+1)\gdom\t$. Let
$\alpha=\frac{(q-1)\res_\t(i)}{\res_\t(i)-\res_\s(i)}$ and
$\beta=\frac{(q-1)\res_\s(i)}{\res_\s(i)-\res_\t(i)}$. Then 
$f_{\s\s}(T_i-\alpha)=f_{\s\t}$ by 
Proposition~\ref{Tf multiplication}. Similarly, by the left hand analogue of
Proposition~\ref{Tf multiplication} (interchanging the roles of $\s$ and $\t$), 
together with (\ref{gamma properties})(ii), 
$(T_i-\beta)f_{\t\t}=(\gamma_\t/\gamma_\s)f_{\s\t}$. Therefore,
$$(T_i-\beta)F_\t=\frac1{\gamma_\t}(T_i-\beta)f_{\t\t}
                 =\frac1{\gamma_\s}f_{\s\t}
		 =\frac1{\gamma_\s}f_{\s\s}(T_i-\alpha)
		 =F_\s(T_i-\alpha).$$
By induction, there exist elements~$\Phi_\s$ and $\Psi_\s$ which
satisfy properties (i)--(iii). Define $\Psi_\t=\Psi_\s(T_i-\beta)$
and $\Phi_\t=\Phi_\s(T_i-\alpha)$; then, by induction and the
last equation, 
$$\Psi_\t F_\t=\Psi_\s(T_i-\beta)F_\t
                  =\Psi_\s F_\s(T_i-\alpha)
		  =F_\tlam\Phi_\s(T_i-\alpha)
		  =F_\tlam\Phi_\t.$$
Hence, (i) holds. Next, again by induction we have
$$\Psi_\t=\Psi_\s(T_i-\beta)
         =\Big(T_{d(\s)}+\Sum_{v<d(\s)} p_{\s v}T_v\Big)(T_i-\beta)
         =T_{d(\t)}+\Sum_{w<d(\t)} p_{\t w}T_w,$$
by standard properties of the Bruhat order since
$d(\t)=d(\s)(i,i+1)>d(\s)$. This proves~(ii). Finally, using
induction once more (and a quick calculation for the second
equality),
\begin{align*} 
\glam\Phi_\t\PS\Psi_\t^*
    &=\glam\Phi_\s\PS(T_i-\alpha)(T_i-\beta)\Psi_\s^* \\
    &=\frac{(q\res_\s(i)-\res_\t(i))(\res_\s(i)-q\res_\t(i))}
            {(\res_\s(i)-\res_\t(i))^2}\glam \Phi_\s\PS\Psi_\s^* \\
    &= \frac{(q\res_\s(i)-\res_\t(i))(\res_\s(i)-q\res_\t(i))}
            {(\res_\s(i)-\res_\t(i))^2}\gamma_\s=\gamma_\t,
\end{align*}
the last equality coming from (\ref{gamma properties})(ii).
This proves (iii) and so completes the proof.
\end{proof}

We are not claiming that the elements $\Phi_\t$ and $\Psi_\t$ are
uniquely determined by the conditions of the Proposition; ostensibly,
these elements depend upon the choice of reduced expression for $d(\t)$.
In what follows we only need to know that elements with these
properties exist.

\begin{Corollary}
Suppose that $\t$ is a standard $\lambda$--tableau.  Then
\begin{enumerate}
\item $F_\tlam=\frac{\glam}{\gamma_\t}\Psi_\t F_\t\Psi_\t^*$; and,
\item $F_\t=\frac{\glam}{\gamma_\t}\Phi_\t^*F_\tlam\Phi_\t$.
\end{enumerate}\label{factors of F_lam}
\end{Corollary}

\begin{proof}
Using parts (iii) and (i) of the Proposition,
respectively, shows that
$$F_\tlam=\frac{\glam}{\gamma_\t}F_\tlam\Phi_\t\Psi_\t^*
         =\frac{\glam}{\gamma_\t}\Psi_\t F_\t\Psi_\t^*;$$
this proves (i). Part (ii) follows from (i) by `conjugating' (i) by
$\Phi_\t$.
\end{proof}

The main reason why we are interested in $\Psi_\t$ and $\Phi_\t$
is the following.

\begin{Proposition}
Suppose that $\s$ and $\t$ are standard
$\lambda$--tableaux. Then
$$f_{\s\t}=\Phi_\s^*f_{\tlam\tlam}\Phi_\t.$$
\end{Proposition}

\begin{proof}
By the definition of $f_{\s\t}$ and Proposition~\ref{Phi-Psi prop}(ii)
we have 
$$f_{\s\t}=F_\s T_{d(\s)}^* m_\lambda T_{d(\t)} F_\t
	  =F_\s \Psi_\s^* m_\lambda \Psi_\t F_\t
	      -\sum_{(v,w)<(d(\s),d(\t))} 
	      p_{\s v}p_{\t w}F_\s T_v^*m_\lambda T_w F_\t.$$
Now if $(v,w)<(d(\s),d(\t))$ then $T_v^*m_\lambda T_w$ belongs
to the span of the $m_{\u\v}$ where $(\u,\v)\gdom(\s,\t)$.
Therefore, by Proposition~\ref{f_st properties}(i), $T_v^*m_\lambda T_w$
belongs to the span of the $f_{\u\v}$ where either $\u$ and $\v$
are standard $\lambda$--tableaux and $(\u,\v)\gdom(\s,\t)$,
or $\Shape(\u)=\Shape(\v)\gdom\lambda$; consequently, $F_\s
T_v^*m_\lambda T_w F_\t=0$ by Proposition~\ref{f_st properties}(iii). Hence,
by Theorem~\ref{idempotents}(i) and Proposition~\ref{Phi-Psi prop}(i),
$$f_{\s\t} =F_\s \Psi_\s^* m_\lambda \Psi_\t F_\t
           =\Phi_\s^* F_\tlam m_\lambda F_\tlam\Phi_\t
           =\Phi_\s^* f_{\tlam\tlam}\Phi_\t$$
as required.
\end{proof}

Applying this result to $f_{\tlam\tllam}$ shows that
$f_{\tlam\tllam}=\Phi_\tlam^*f_{\tlam\tlam}\Phi_\tllam
=\glam\Flam\Phi_\tllam$, the last equality following because
$\Phi_{\tlam}=1$. Using Proposition~\ref{Phi-Psi prop}(iii) to multiply this
equation on the right by $\Psi_{\tllam}^*$, and recalling
Corollary~\ref{z=f}, now yields the following.

\begin{Corollary}
Suppose that $\lambda$ is a multipartition of $n$. Then
$$\Flam=\frac1{\gamma_\tllam} f_{\tlam\tllam}\Psi_{\tllam}^*
  =\frac1{\gamma_\tllam\gamma_\tlamp'}z_\lambda\Psi_{\tllam}^*.$$
\end{Corollary}

By Lemma~\ref{generic degree lemma} in order to compute the Schur elements
it suffices to calculate $\tau(\Flam)$ for each $\lambda$; so we are
reduced to finding $\tau(z_\lambda\Psi_{\tllam}^*)$. To do this we
rewrite $z_\lambda$ with respect to the Ariki--Koike basis.  In types
$A$ and $B$ (that is, $r=1$ or $r=2$) this is reasonably
straightforward; in general we have to work much harder.

Until further notice, fix a multipartition $\lambda$ and let
$\barl=\rtuple\barl$ be the multipartition with $\barl
^{(s)}=(|\lambda^{(s)}|)$ for $1\le s\le r$. Then $\wbarl=d(\t^\barl)$
is the permutation
$$\Big(\begin{array}{ccc|ccc|c|ccc}
 1&\cdots&\a_2&\a_2+1&\cdots &\a_3&\cdots& \a_r+1&\cdots&n\\
 \b_1+1&\cdots&n&\b_2+1& \cdots&\b_1&\cdots&1&\cdots&\b_{r-1}
\end{array}\Big),$$
where $\a_s=\sum_{t=1}^{s-1}|\lambda^{(t)}|$ and
$\b_s=\sum_{t=s+1}^r|\lambda^{(t)}|$ for $1\le s\le r$.

Now, $\t^\barl\wbarl$ is a standard tableau; therefore, $\wbarl$ is a
distinguished right coset representative for~$\Sym_\barl$ in~$\Sym_n$
(for example, by \cite[Prop.~3.3]{M:ULect}). Equivalently, $\wbarl$ is
a  distinguished left coset representative for $\Sym_{\overline\lamp}$
in $\Sym_n$ since $\Sym_\barl w_\barl=w_\barl\Sym_\lamp$.
Consequently, there exists $w_{\lambda/\barl}\in\Sym_{\overline\lamp}$
such that $w_\lambda=\wbarl w_{\lambda/\barl}$ and
$\len(w_\lambda)=\len(\wbarl)+\len(w_{\lambda/\barl})$. To proceed we
need to factorize~$\wbarl $. 

Suppose that $a$ and $b$ are non--negative integers let
$w_{a,b}=w_{((a),(b))}$; so, $\t_{((a),(b))}=\t^{((a),(b))}w_{a,b}$.
More concretely, $w_{a,0}=1=w_{0,b}$ and if $a>0$ and $b>0$ then
$$w_{a,b}=\Big(\begin{array}{*4c|*4c} 1&2&\cdots&a&a+1&a+2&\cdots&a+b\\
	                   b+1&b+2&\cdots&a+b&1&2&\cdots&b
	         \end{array}\Big).$$
If $i\le j$ we also set $s_{i,j}=s_is_{i+1}\dots s_j$ and
$s_{j,i}=1$; so $s_{i,j}=(j+1,j,\dots,i)$. For
convenience let $T_{i,j}=T_{s_{i,j}}$. It is not hard to see that
$w_{a,b}=(s_{a+b-1,1})^b$.

The permutations $w_{a,b}$ were studied by Dipper and
James~\cite{DJ:B}; in particular, they observed 
that $w_{a,b}=s_{a,a+b-1}w_{a-1,b}$, with the lengths adding, which
implies the following. (Note that our $s_{i,j}$ is Dipper
and James' $s_{i,j-1}$ when~$i<j$; also if $i\ge j$ then Dipper and
James set $s_{i,j}=(s_{j,i-1})^{-1}$ whereas we have $s_{i,j}=1$.)
We let $\Sym_{(b,a)}=\Sym_b\times\Sym_a\incl\Sym_n$ (natural embedding).

\begin{Lemma}
Suppose that $a$ and $b$ are positive integers with $a+b\le n$. 
Then $w_{a,b}=s_{a,a+b-1}\dots s_{1,b}$ and
$\len(w_{a,b})=\len(s_{a,a+b-1})+\dots+\len(s_{1,b})$.
Moreover, $w_{a,b}$ is a distinguished left coset representative for
$\Sym_{(b,a)}$ in $\Sym_n$; that is,
$\len(w_{a,b}v)=\len(w_{a,b})+\len(v)$ for all
$v\in\Sym_{(b,a)}$.
\label{w_a reduced}
\end{Lemma}

For $s=1,\dots,r$ define $w_{\barl,s}=w_{n_s,\b_s}$, where
$n_s=|\lambda^{(s)}|$ and $\b_s=\sum_{t=s+1}^r|\lambda^{(t)}|$; in
particular, $\b_0=n$ and $w_{\barl,r}=1$.  Now $n_s+\b_s=\b_{s-1}$ so
$w_{\barl,s}\in\Sym_{\b_{s-1}}$ for each~$s$; therefore,
$w_{\barl,s}\dots w_{\barl,r-1}$ is an element of $\Sym_{\b_{s-1}}$
and, consequently,
$$\len(w_{\barl,1}\dots w_{\barl,r-1})
        =\len(w_{\barl,1})+\len(w_{\barl,2}\dots w_{\barl,r-1})
        =\dots=\len(w_{\barl,1})+\dots+\len(w_{\barl,r-1})$$
by Lemma~\ref{w_a reduced}.  Noting that 
$(k)w_\barl=(k)w_{\barl,1}\dots w_{\barl,s}$
for $k=1,2,\dots,\a_{s+1}$ we have shown the following.

\begin{Point}
{\it}* Let $\lambda$ be a partition of $n$. Then
$w_\lambda=w_\barl w_{\lambda/\barl}
           =w_{\barl,1}\dots w_{\barl,r-1}w_{\lambda/\barl}$.
Moreover, 
$\len(w_\lambda)=\len(w_\barl)+\len(w_{\lambda/\barl})
                =\len(w_{\barl,1})+\dots+\len(w_{\barl,r-1})
                  +\len(w_{\lambda/\barl}).$
\label{w_lam factorization}\end{Point}

The point of the lengths adding is that 
$T_{w_\lambda}=T_{w_{\barl,1}}\dots T_{w_{\barl,r-1}}
                   T_{w_{\lambda/\barl}}$.
We will use this below without further comment.

Although we won't need it notice that $\len(w_{a,n-a})=a(n-a)$
by~Lemma~\ref{w_a reduced} (or directly); so (\ref{w_lam factorization})
implies that 
$\len(w_\barl)=\sum_{1\le s<t\le r}|\lambda^{(s)}||\lambda^{(t)}|$.

\begin{Example}
Let $\lambda=\((2,1^2),(2,1),(2)\)$. Then
$\barl=\((4),(3),(2)\)$ and
\begin{align*}
w_\lambda&=\Big(\begin{array}{*4r|*3r|*2r} 1& 2& 3& 4& 5& 6& 7& 8& 9\\
                                   6&9&7& 8& 3& 5& 4& 1& 2
	           \end{array}\Big)=(1,6,5,3,7,4,8)(2,9),\\
\wbarl&=\Big(\begin{array}{*4r|*3r|*2r} 1& 2& 3& 4& 5& 6& 7& 8& 9\\
                              6& 7& 8& 9& 3& 4& 5& 1& 2
              \end{array}\Big)=(1,6,4,9,2,7,5,3,8);
\end{align*}
therefore, $w_{\lambda/\barl}=(7,9,8)(4,5)$. Further,
$\b=(5,2,0)$ so, using~Lemma~\ref{w_a reduced},
\begin{align*}
w_{\barl,1}&=w_{4,5}
            =\Big(\begin{array}{*4r|*5r} 1& 2& 3& 4& 5& 6& 7& 8& 9\\
                               6& 7& 8& 9& 1& 2& 3& 4& 5
              \end{array}\Big)\\
           &=(1,6,2,7,3,8,4,9,5)=s_{4,8}s_{3,7}s_{2,6}s_{1,5}\\
w_{\barl,2}&=w_{3,2}=\Big(\begin{array}{*3r|*2r} 1& 2& 3& 4& 5\\
                                       3& 4& 5& 1& 2
              \end{array}\Big)=(1,3,5,2,4)=s_{3,4}s_{2,3}s_{1,2}
\end{align*}
and $w_{\barl,3}=1$. It is evident that 
$w_\lambda=w_{\barl,1}w_{\barl,2}w_{\lambda/\barl}$.
\end{Example}

Given $i\le j$ and $Q\in R$ let $L_{i,j}(Q)=(L_i-Q)\dots(L_j-Q)$; if
$i>j$ set $L_{i,j}(Q)=1$.  In particular,
$u^-_\lamp=L_{1,\b_1}(Q_1)\dots L_{1,\b_{r-1}}(Q_{r-1})$ where
$\b_s=|\lambda^{(s+1)}|+\dots+|\lambda^{(r)}|$, as above.

As a final piece of notation, given $0\le l\le m<n$ let $H_{k,m}$ be
the $R$--submodule of~$\H$ spanned by the elements
$\set{L_1^{c_1}\dots L_k^{c_k}T_w|0\le c_i<r\And w\in\Sym_m}$. In
general~$H_{k,m}$ is neither a subalgebra nor a submodule of~$\H$;
however, $H_{k,m}$ is a right $\H(\Sym_m)$--module. 

One of the difficulties in working with $\H$ is that, in general, the
powers of the $L_k$ are not necessarily linear combinations of
products of the $L_i^s$, for $1\le i\le k$ and $1\le s<r$; however, it
is always true that $L_k^s\in H_{k,m}$ for all $s\ge0$ provided that
$m\ge k$ (since in this case $\<T_0,\dots,T_{k-1}\>\ss H_{k,m}$ 
by~(\ref{AK basis})).

\begin{Lemma}
\label{H_km}
Suppose that $h\in H_{k-1,k}$ for some $k$ with $1\le k\le n$. 
\begin{enumerate}
\item If $1\le i<n$ then
$$T_ih\in\begin{cases} H_{k-1,k},&\If 1\le i<k-1,\\
                        H_{k,k},&\If i=k-1,\\
                        H_{k-1,i+1},&\If k\le i<n.
\end{cases}$$
\item If $1\le s<r$ then $hL_k^s=\sum_{t=0}^sL_k^th_t$ for some 
$h_t\in H_{k-1,k}$.
\end{enumerate}
\end{Lemma}

\begin{proof}
By definition $h$ is a linear combination of terms of the form
$L_{k-1}^{c_{k-1}}\dots L_1^{c_1}T_w$, for some $0\le c_i<r$ and
$w\in \Sym_k$. By (\ref{L-comm})(ii), if $i\ne m-1,m$ then
$T_iL_m=L_mT_i$; whereas $T_mL_m=L_{m+1}(T_m-q+1)$ and
$T_{m-1}L_m=L_{m-1}T_{m-1}+(q-1)L_m$. Combining these facts proves (i).

For part (ii), observe that if $w\in\Sym_k$ then $T_wL_k=L_kT_w$
unless $s_{k-1}$ appears in a reduced expression of $w$; however,
$T_{k-1}L_k=L_{k-1}T_{k-1}+(q-1)L_k$. Part~(ii) now follows by
induction on~$s$ using the remarks of the last paragraph.
\end{proof}

\begin{Lemma}
Fix $Q\in R$ and let $i$, $j$, $c$ and $d$ be positive
integers with $i\le c\le d\le j<n$. Then 
$$T_{i,j}L_{c,d}(Q)=\Big\{L_{c+1,d+1}(Q)
  +\sum_{k=c+1}^{d+1}L_kh_k\Big\}T_{i,j}$$
for some elements $h_k\in H_{k-1,k}$ for $c+1\le k\le d+1$.
\label{T_ijL_cd}
\end{Lemma}

\begin{proof}
If $d=j$ then by convention $T_{d+1,j}=1$. Therefore, by 
(\ref{L-comm})(iv) we have
\begin{align*}
T_{i,j}L_{c,d}(Q)
  &=T_{i,c-1}T_{c,d}(L_c-Q)\dots(L_d-Q)T_{d+1,j}\\
  &=T_{i,c-1} T_c(L_c-Q)\dots T_d(L_d-Q)T_{d+1,j}.\\
\intertext{Now, $T_kL_k=qL_{k+1}T_k^{-1}=L_{k+1}T_k-(q-1)L_{k+1}$, so}
T_{i,j}L_{c,d}(Q)
  &=T_{i,c-1}\prod_{k=c+1}^{d+1}\Big\{(L_k-Q)T_{k-1}-(q-1)L_k\Big\}
       \cdot T_{d+1,j},
\intertext{where we read the terms from left to right with increasing
values of $k$ --- the order of the factors is important here. Now
$T_{k-1}$ commutes with $L_m$ if $k<m$; therefore, each $T_{k-1}$ in
this product commutes with all of $L_m$'s which appear to its right.
A straightforward induction on $(d-c)$ using (\ref{L-comm}) shows that 
there exist $h_k'\in H_{k-1,k}$ such that}
T_{i,j}L_{c,d}(Q)
   &=T_{i,c-1}\Big\{L_{c+1,d+1}(Q)T_{c,d}
     +\sum_{k=c+1}^{d+1}L_kh_k'T_{k,d}\Big\}T_{d+1,j}.
\end{align*}
(Recall that our convention is that $T_{d+1,d}=1$.)
If~$c+1\le k\le d+1$ then by (\ref{L-comm})
$$T_{i,c-1}L_kh_k'T_{k,d}T_{d+1,j}
   =L_kT_{i,c-1}h_k'T_{k,j}=L_kT_{i,c-1}h_k'T_{i,k-1}^{-1}T_{i,j}.$$ 
Also, using~(\ref{L-comm})(ii), 
$T_{i,c-1}L_{c+1,d+1}(Q)=L_{c+1,d+1}(Q)T_{i,c-1}$. Therefore,
$$T_{i,j}L_{c,d}(Q)
   =L_{c+1,d+1}(Q)T_{i,j}+\sum_{k=c+1}^{d+1}L_kh_kT_{i,j},
$$
where $h_k=T_{i,c-1}h_k'T_{i,k-1}^{-1}\in H_{k-1,k}$ by Lemma~\ref{H_km},
for $c+1\le k\le d+1$.
\end{proof}

This brings us to the key technical lemma.

\begin{Lemma}
Fix $Q\in R$ and let $a$, $b$, $c$ and $d$ be non--negative
integers such that $a+b\le n$ and $1\le c\le d\le b$. Then 
there exist elements $h_k\in H_{k-1,k}$ such
that
$$T_{w_{a,b}}L_{c,d}(Q)
  =\Big\{L_{a+c,a+d}(Q)+\sum_{k=a+c}^{a+d}L_kh_k\Big\} T_{w_{a,b}}.$$
\label{technical lemma}
\end{Lemma}

\begin{proof}
We argue by induction on $a$. If $a=0$ then $w_{a,b}=1$ and
there is nothing to prove; so suppose that $a\ge1$. Then
$w_{a,b}=s_{a,a+b-1}w_{a-1,b}$ by Lemma~\ref{w_a reduced}. Therefore, by
induction and Lemma~\ref{T_ijL_cd}, respectively, there exist elements 
$h_k',h_k''\in H_{k-1,k}$ such that
\begin{align*}
T_{w_{a,b}}L_{c,d}(Q)
 &=T_{a,a+b-1}T_{w_{a-1,b}}L_{c,d}(Q)\\
 &=T_{a,a+b-1}\Big\{L_{a+c-1,a+d-1}(Q)
  +\sum_{k=a+c-1}^{a+d-1}L_kh_k'\Big\}T_{w_{a-1,b}}\\
 &=\Big\{L_{a+c,a+d}(Q)
     +\sum_{k=a+c}^{a+d}L_kh_k''\Big\}T_{a,a+b-1}T_{w_{a-1,b}}\\
 &\qquad\qquad+\sum_{k=a+c-1}^{a+d-1}T_{a,a+b-1}L_kh_k'T_{w_{a-1,b}}.
\end{align*}
If $a+c-1\le k\le a+d-1$ then by (\ref{L-comm})
\begin{align*}
T_{a,a+b-1}L_kh_k'T_{w_{a-1,b}}
  &=T_{a,k-1}T_kL_kT_{k+1,a+b-1}h_k'T_{w_{a-1,b}}\\
  &=qT_{a,k-1}L_{k+1}T_k^{-1}T_{k+1,a+b-1}h_k'T_{w_{a-1,b}}\\
  &=qL_{k+1}T_{a,k-1}T_k^{-1}h_k'T_{k+1,a+b-1}T_{w_{a-1,b}}\\
  &=qL_{k+1}T_{a,k-1}T_k^{-1}h_k'T_{a,k}^{-1}T_{a,a+b-1}T_{w_{a-1,b}}\\
  &=L_{k+1}T_{a,k-1}(T_k-q+1)h_k'T_{a,k}^{-1}T_{w_{a,b}},
\end{align*}
where the third equality follows by (\ref{L-comm}) because
$h_k'\in H_{k-1,k}$. By Lemma~\ref{H_km} the element
$T_{a,k-1}(T_k-q+1)h_k'T_{a,k}^{-1}$ 
belongs to $H_{k,k+1}$; so, combining these equations proves the Lemma.
\end{proof}

Lemma~\ref{technical lemma} is enough to compute the Schur elements when
$r=2$; in order to cover the general case we delicately apply
Lemma~\ref{technical lemma} several times.

\begin{Corollary}
Suppose that $a$, $b$, $c$ and $s$ are non--negative integers with
$a+b\le n$, $1\le c\le b$ and $1\le s<r$. Then 
$T_{w_{a,b}}L_c^s=\sum_{t=1}^sL_{a+c}^th_tT_{w_{a,b}}$
some $h_t\in H_{a+c-1,a+c}$.
\label{T_wabL_c}
\end{Corollary}

\begin{proof}
If $s=1$ this is the special case of Lemma~\ref{technical lemma}
corresponding to the choices $d=c$ and $Q=0$; in particular, this
implies that the exponent $t$ is always at least~$1$. If $s>1$ then
the result follows by induction using Lemma~\ref{H_km}(ii). 
\end{proof}

For each integer $k$ let $c_\lambda(k)=s$ if $k$ appears in
component~$s$ of $\tlam$. Observe that for all~$k$ we have
$c_\lambda(k)=\min\set{1\le s\le r|
           k\le |\lambda^{(1)}|+\dots+|\lambda^{(s)}|}$.
For the next Lemma recall that 
$\a_s=|\lambda^{(1)}|+\dots+|\lambda^{(s-1)}|$ and
$\b_s=|\lambda^{(s+1)}|+\dots+|\lambda^{(r)}|$
for $1\le s\le r$.

\begin{Lemma}
Suppose that $\lambda$ is a multipartition of $n$. Then
$$T_{w_\lambda}u^-_\lamp=
  L_{\a_2+1,n}(Q_1)\dots L_{\a_r+1,n}(Q_{r-1}) T_{w_\lambda}
     +\sum_{k=1}^n\sum_{c=1}^{c_\lambda(k)-1}L_k^ch_{kc}$$
for some $h_{kc}\in H_{k-1,n}$.
\label{rewriting}
\end{Lemma}

\begin{proof}
First note that if $r=1$ then $u_\lamp^-=1$ and there is nothing 
to prove; so we assume that~$r>1$. Next,
$T_{w_\lambda}u^-_\lamp=T_\wbarl T_{w_{\lambda/\barl}}u^-_\lamp
                      =T_\wbarl u^-_\lamp T_{w_{\lambda/\barl}}$
by part~(iv) of (\ref{L-comm}) since
$w_{\lambda/\barl}\in\Sym_{\overline{\lamp}}$.
Consequently, it is enough to consider $T_\wbarl u_\lamp^-$; 
equivalently, we may assume that $\lambda=\barl$ (for clarity we will
continue to write $\barl$). 

By (\ref{w_lam factorization}), $w_\barl=w_{\barl,1}\dots w_{\barl,r-1}$ 
where $w_{\barl,s}\in\Sym_{\b_{s-1}}$. Therefore, 
by (\ref{L-comm})(iv),
\begin{align*}
  T_{w_\barl}u^-_\lamp
   &=T_{w_{\barl,1}}\dots T_{w_{\barl,r-1}}
            L_{1,\b_1}(Q_1)\dots L_{1,\b_{r-1}}(Q_{r-1})\\
   &=T_{w_{\barl,1}}L_{1,\b_1}(Q_1)\dots 
            T_{w_{\barl,r-1}}L_{1,\b_{r-1}}(Q_{r-1}).
\end{align*}
We also let $n_s=|\lambda^{(s)}|$.
Suppose that $1\le s\le r-1$. For $k=1,\dots,\b_{s-1}$ set
$$c_{\lambda,s}(k)=\min\set{0\le t\le r-s|
   k\le|\lambda^{(s)}|+\dots+|\lambda^{(s+t)}|}$$
and let $H^{(s)}$ be the $R$--submodule of $\H$ spanned by
elements of the form $L_k^ch_{kc}$ where $n_s<k\le\b_{s-1}$, 
$1\le c\le c_{\lambda,s}(k)$ and $h_{kc}\in H_{k-1,\b_{s-1}}$. 
We claim that
$$T_{w_{\barl,s}}L_{1,\b_s}(Q_s)\dots 
            T_{w_{\barl,r-1}}L_{1,\b_{r-1}}(Q_{r-1})
      \equiv P^{(s)}\pmod{H^{(s)}},$$
where $P^{(s)}=\prod_{t=s}^{r-1}L_{n_s+\dots+n_t+1, \b_{s-1}}(Q_t)\cdot
            T_{w_{\barl,s}}\dots T_{w_{\barl,r-1}}$. 
Taking $s=1$ will prove the Lemma because $\b_0=n$ and
$c_\lambda(k)=c_{\lambda,1}(k)+1$ by the remarks before the Lemma.
Note that $P^{(s)}$ has degree $c_{\lambda,s}(k)$ when we consider it
as a polynomial in~$L_k$.

To prove the claim we argue by downwards induction on $s$. If $s\ge
c_\lambda(n)$ then $\b_s=0$ so that
$T_{w_{\barl,s}}=L_{1,\b_s}(Q_s)=P^{(s)}=1$ and there is nothing to
prove. Suppose then that $s<c_\lambda(n)$. Then, by induction,
$$T_{w_{\barl,s}}L_{1,\b_s}(Q_s)\dots 
            T_{w_{\barl,r-1}}L_{1,\b_{r-1}}(Q_{r-1})
  = T_{w_{\barl,s}}L_{1,\b_s}(Q_s)\Big\{P^{(s+1)}+h^{(s+1)}\Big\}
$$
for some $h^{(s+1)}\in H^{(s+1)}$. Now $\b_s>0$ since $s<c_\lambda(n)$, 
so $\b_{s-1}=\b_s+n_s\ge\b_s\ge1$. Hence, by Lemma~\ref{technical lemma},
$$\begin{array}{r@{}l}
T_{w_{\barl,s}}L_{1,\b_s}(Q_s)&\Big\{P^{(s+1)}+h^{(s+1)}\Big\}
   = T_{w_{n_s,\b_s}}L_{1,\b_s}(Q_s)\Big\{P^{(s+1)}+h^{(s+1)}\Big\}\\
   =\,&\Big\{L_{n_s+1,\b_{s-1}}(Q_s)T_{w_{n_s,\b_s}}
       +\!\!\Sum_{k=n_s+1}^{\b_{s-1}}L_kh_kT_{w_{n_s,\b_s}}\Big\}
      \Big\{P^{(s+1)}+h^{(s+1)}\Big\},
\end{array}$$
for some $h_k\in H_{k-1,k}$. We will move $T_{w_{\barl,s}}$
past $P^{(s+1)}$ and $h^{(s+1)}$ using Lemma~\ref{technical lemma}.

Viewing $P^{(s+1)}$ as a polynomial in the $L_m$'s, if
$n_{s+1}<m\le\b_s$ then the degree of $L_m$ in~$P^{(s+1)}$ is at most
$c_{\lambda,s+1}(m)$; further, $L_m$ does not appear in $P^{(s+1)}$ if
either~$m\le n_{s+1}$ or $m>\b_s$.  Let
$p^{(s+1)}=\prod_{t=s+1}^{r-1}L_{n_{s+1}+\dots+n_t+1,\b_s}(Q_t)$. Then
by Lemma~\ref{technical lemma} and Corollary~\ref{T_wabL_c} there exist $h_{k,1},
h_{ck,2}\in H_{k-1,k}$ such that
\begin{align*}
T_{w_{n_s,\b_s}}p^{(s+1)}
   &=\Big\{L_{n_s+n_{s+1}+1,\b_{s-1}}(Q_{s+1})
         +\Sum_{k=n_s+n_{s+1}+1}^{b_{s-1}}L_kh_{k,1}\Big\} 
           T_{w_{n_s,b_s}}p^{(s+2)}\\
   &=\Big\{L_{n_s+n_{s+1}+1,\b_{s-1}}(Q_{s+1})
           L_{n_s+n_{s+1}+n_{s+2}+1,\b_{s-1}}(Q_{s+2})\\
   &\hspace*{12mm}+\sum_{k=n_s+n_{s+1}+1}^{\b_{s-1}}
          \sum_{c=1}^{c_2(k)}L_k^ch_{ck,2}
          \Big\}T_{w_{n_s,b_s}}p^{(s+3)},\\
\end{align*}
where $c_2(k)=1$ if $n_s+n_{s+1}<k\le n_s+n_{s+1}+n_{s+2}$ and
$c_2(k)=2$ if $n_s+n_{s+1}+n_{s+2}<k\le\b_{s-1}$. For the second
equality we have used Lemma~\ref{technical lemma}, Lemma~\ref{H_km} and
Corollary~\ref{T_wabL_c}.  Looking at the definition of $c_{\lambda,s}(m)$, if
$n_{s+1}<m\le\b_s$ then $c_{\lambda,s+1}(m)=c_{\lambda,s}(m+n_s)-1$;
therefore, continuing in this way we see that
$$T_{w_{n_s,\b_s}}p^{(s+1)}
   \equiv\prod_{t=s+1}^{r-1}
          L_{n_s+n_{s+1}+\dots+n_t+1,\b_{s-1}}(Q_t)\cdot T_{w_{n_s,\b_s}}
      \pmod{ H^{(s)} }.$$
In obtaining this equation notice that if $m>k$ then
$L_kh_kL_mh_k=L_m(L_kh_kh_m)$ and $L_kh_kh_m\in H_{m-1,m}$, for any
$h_k\in H_{k-1,k}$ and $h_m\in H_{m-1,m}$; similarly, if $m<k$ then
$h_kL_mh_m\in H_{k-1,k}$ by Lemma~\ref{H_km}.  It follows that if $L_k^ch$
appears in this expansion, for some $h\in H_{k-1,k}$, then
$n_s+n_{s+1}<k\le\b_{s-1}$ and so~$1\le c\le
c_{\lambda,s+1}(k-n_s)=c_{\lambda,s}(k)-1$ by the remarks above;
hence, $L_kh_kL_mh_m\in H^{(s)}$ for all $k$ and $m$.

Now consider a term from the inductive step of the form 
$L_kh_kT_{w_{n_s,\b_s}}P^{(s+1)}$, where $n_s<k\le\b_{s-1}$ and
$h_k\in H_{k-1,k}$. What we have just shown combined with 
Lemma~\ref{H_km}(ii) shows that $L_kh_kT_{w_{n_s,\b_s}}P^{(s+1)}$ is equal
to a linear combination of terms of the form $L_m^ch_{md}$,
where $n_s<m\le\b_{s-1}$, $1\le c\le c_{\lambda,s}(m)$ and 
$h_{md}\in H_{m-1,m}$. Moreover, $c\le c_{\lambda,s}(m)$ with equality
only if $k=m>n_s+n_{s+1}$ and $c=1=c_{\lambda,s}(m)$ if 
$n_s<m\le n_s+n_{s+1}$. Hence, 
$L_kh_kT_{w_{n_s,\b_s}}P^{(s+1)}\in H^{(s)}$.

Therefore, combining the last two paragraphs we have shown that
$$T_{w_{n_s,\b_s}}L_{1,\b_s}(Q_s)P^{(s+1)}
  \equiv P^{(s)}\pmod{H^{(s)}}.$$
By similar arguments, the terms
$L_{n_s+1,\b_{s-1}}(Q_s)T_{w_{n_s,\b_s}}h^{(s+1)}$ and
$L_kh_kT_{w_{n_s,\b_s}}h^{(s+1)}$ from the inductive step also
belong to~$H^{(s)}$. We leave the details to the reader.
\end{proof}

We are now basically done. The next result essentially computes
$\tau(z_\lambda \Psi_\tllam^*)$; we record it separately because it is
also the key to showing that $M^\lambda$ is a self--dual $\H$--module
and that $\tilde S^{\lambda'}$ is isomorphic to the
dual of $S^\lambda$; see \cite{M:tilting}.

\begin{Proposition}
Suppose that $\lambda$ is a multipartition of $n$. Then
$$\tau(z_\lambda T_{w_{\lambda}}^*)
 =(-1)^{n(r-1)}q^{\len(w_\lambda)}\prod_{s=1}^rQ_s^{n-|\lambda^{(s)}|}.$$
\label{tau(z_lambda)}
\end{Proposition}

\begin{proof}
Let 
$\tilde u_\lamp^-=L_{\a_2+1,n}(Q_1)\dots L_{\a_r+1,n}(Q_{r-1})$. 
Using (\ref{tau properties})(i) and Lemma~\ref{rewriting} shows that
\begin{align*}
\tau(z_\lambda T_{w_{\lambda}}^*)
    &=\tau\(m_\lambda T_{w_\lambda}n_\lamp T_{w_{\lambda}}^*\)
     =\tau\(x_\lambda u_\lambda^+ T_{w_\lambda}u_\lambda^-
              y_\lamp T_{w_{\lambda}}^*\)\\
    &=\tau\(u_\lambda^+ T_{w_\lambda}u_\lambda^-
              y_\lamp T_{w_{\lambda}}^*x_\lambda\)\\
    &=\tau\(u_\lambda^+\tilde u_\lamp^-
             T_{w_\lambda}y_\lamp T_{w_{\lambda}}^*x_\lambda\)
        +\tau\(u_\lambda^+hy_\lamp T_{w_{\lambda}}^*x_\lambda\),
\end{align*}
where $h=\sum_{k=1}^n\sum_{c=1}^{c_\lambda(k)-1}L_k^ch_{kc}$ for some
$h_{kc}\in H_{k-1,n}$. Consider $u_\lambda^+$ as a polynomial in
$L_k$. Then the degree of~$L_k$ in $u_\lambda^+$ is $r-c_\lambda(k)$.
Therefore, if $h_{kc}\ne0$ then $u_\lambda^+L_k^ch_{kc}$ is a
polynomial in~$L_k$ which is left divisible by $L_k$ and has degree at
most $r-1$. If~$m>k$ then $L_m$ appears
in~$u_\lambda^+L_k^ch_{kc}$ with exponent at most $r-c_\lambda(m)<r$.
(If~$m<k$ then $L_m^d$ can appear in $u_\lambda^+L_k^ch_{kc}$ for
$d\ge r$; however, this does not matter because such terms can be written
as a linear combination of Ariki--Koike basis elements
in~$H_{k-1,k}$.) Therefore, 
$\tau\(u_\lambda^+L_k^ch_{kc}y_\lamp T_{w_{\lambda}}^*x_\lambda\)=0$ 
by (\ref{tau properties})(ii) and, consequently,
$\tau\(u_\lambda^+hy_\lamp T_{w_{\lambda}}^*x_\lambda\)=0$; hence,
$$\tau(z_\lambda T_{w_{\lambda}}^*)
  =\tau\(u_\lambda^+\tilde u_\lamp^-
           T_{w_\lambda}y_\lamp T_{w_{\lambda}}^*x_\lambda\).$$
Considered as polynomials in $L_k$,
$u_\lambda^+$ has degree $r-c_\lambda(k)$ and $\tilde u_\lamp^-$ has
degree $c_\lambda(k)-1$; consequently, each $L_k$ has degree $r-1$ in
$u_\lambda^+\tilde u_\lamp^-$. Therefore, by (\ref{tau properties})
again,
\begin{align*}
 \tau(z_\lambda T_{w_{\lambda}}^*)
      &=\tau(u_\lambda^+\tilde u_\lamp^-)\cdot
        \tau(T_{w_\lambda}y_\lamp T_{w_{\lambda}}^*x_\lambda)\\
      &=\prod_{s=2}^r(-Q_s)^{\a_s}\cdot
        \prod_{s=1}^{r-1}(-Q_s)^{\b_s}\cdot
        \tau(T_{w_\lambda}y_\lamp T_{w_{\lambda}}^*x_\lambda)\\
      &=(-1)^{n(r-1)}\prod_{s=1}^rQ_s^{n-|\lambda^{(s)}|}\cdot
        \tau(T_{w_\lambda}y_\lamp T_{w_{\lambda}}^*x_\lambda),\\
\end{align*}
since $\a_s+\b_s=n-|\lambda^{(s)}|$ for all $s$ (and $\a_1=0=\b_r$).
To complete the proof recall that
$\Sym_\lambda\cap{}^{w_\lambda}\Sym_\lamp=\{1\}$ and that ~$w_\lambda$
is a distinguished $(\Sym_\lambda,\Sym_\lamp)$--double coset
representative (consider the tableaux $\tlam$ and $\tllam$), so 
$$T_{w_\lambda}y_\lamp T_{w_{\lambda}}^*x_\lambda 
     =\sum_{\substack{u\in\Sym_\lamp\\v\in\Sym_\lambda}}
            (-q)^{-\len(u)}T_{w_\lambda}T_uT_{w_\lambda}^*T_v
     =\sum_{\substack{u\in\Sym_\lamp\\v\in\Sym_\lambda}}
            (-q)^{-\len(u)}T_{w_\lambda}T_{uw_\lambda^{-1} v}.
$$
Therefore,
$\tau(T_{w_\lambda}y_\lamp T_{w_{\lambda}}^*x_\lambda)
        =q^{\len(w_\lambda)}$
by (\ref{tau properties})(ii)  (corresponding to $u=v=1$), so the 
Proposition follows.
\end{proof}

We can now give our first formula for the Schur elements
$s_\lambda(\Q)$ of $\H$.

\begin{Corollary}
Suppose that $\lambda$ is a multipartition of $n$. Then 
$$s_\lambda(\Q)=(-1)^{n(r-1)}
  q^{-\len(w_\lambda)}{\gamma_\tllam\gamma_\tlamp'}
   \prod_{s=1}^rQ_s^{|\lambda^{(s)}|-n}.$$
\label{Schur elements}
\end{Corollary}

\begin{proof}
By Lemma~\ref{generic degree lemma},
$s_\lambda(\Q)=1/{\tau(F_\tllam)}$ and
$F_\tlam=\frac1{\gamma_\tllam\gamma_\tlamp'}z_\lambda\Psi_{\tllam}^*$
by Corollary~\ref{factors of F_lam}; therefore,
$s_\lambda(\Q)
      ={\gamma_\tllam\gamma_\tlamp'}/\tau(z_\lambda\Psi_\tllam^*)$.
Now, $\tau$ is a trace form by (\ref{tau properties})(i), so
$$\tau(z_\lambda\Psi_\tllam^*)
  =\tau(x_\lambda u_\lambda^+ T_{w_\lambda}u_\lamp^- 
             y_\lamp\Psi^*_\tllam)
  =\tau(u_\lambda^+T_{w_\lambda}u_\lamp^- 
             y_\lamp\Psi^*_\tllam x_\lambda).
$$
It is well--known and easy to check that the permutation~$w_\lambda$
has the ``trivial intersection property''; that is,
$\Sym_\lambda\cap{}^w\Sym_\lamp\ne\{1\}$ if and only if
$\Sym_\lambda w\Sym_\lamp=\Sym_\lambda w_\lambda\Sym_\lamp$. Therefore,
$y_\lamp T_w^* x_\lambda\ne0$ if and only if
$w\in\Sym_\lambda w_\lambda\Sym_\lamp$ (see, for example, 
\cite[(4.9)]{DJ:reps}). Now,
$\Psi^*_\tllam=T_{w_\lambda}^*+\sum_{w<w_\lambda}p_{\tllam w}T_w^*$ 
by Proposition~\ref{Phi-Psi prop}(ii); so
$y_\lamp\Psi_\tllam^* x_\lambda=y_\lamp T_{w_\lambda}^* x_\lambda$
since $w_\lambda$ is the unique element of minimal length in 
$\Sym_\lambda w_\lambda\Sym_\lamp$. Therefore, 
$$\tau(z_\lambda\Psi_\tllam^*)
  =\tau(u_\lambda^+T_{w_\lambda}u_\lamp^- 
             y_\lamp T_{w_\lambda}^*x_\lamp)
  =\tau(x_\lambda u_\lambda^+T_{w_\lambda}u_\lamp^- 
             y_\lamp T_{w_\lambda}^*)
  =\tau(z_\lambda T_{w_\lambda}^*).$$ 
The result now follows from Proposition~\ref{tau(z_lambda)}.
\end{proof}

A closed formula for $\gamma_\tlamp'$ is given by 
(\ref{gamma properties})(i); therefore, in order to find an explicit
formula for $s_\lambda(\Q)$ we need only compute $\gamma_\tllam$.
Although the formula below looks formidable, its proof follows readily
enough from the definition of~$\gamma_\tllam$.

Recall that the $ij$th {\sf hook} in the diagram $[\lambda^{(s)}]$ is the
collection of nodes to the right of and below the node $(i,j,s)$,
including the node $(i,j,s)$ itself. The $ij$th {\sf hook length}
$h^{\lambda^{(s)}}_{ij}=\lambda^{(s)}_i+\lambda^{(s)'}_j-i-j+1$ is the
number of nodes in the $ij$th hook and the {\sf leg length},
$\ell^{\lambda^{(s)}}_{ij}=\lambda^{(s)'}_j-j+1$, is the number of nodes
in the ``leg'' of this hook. Observe that if $(a,b,c)$ and $(i,j,c)$
are two removable nodes in $[\lambda^{(c)}]$ with $a\le i$ 
and~$j\le b$ then $h^{\lambda^{(c)}}_{aj}=b-a-j+i+1$.

\begin{Lemma}
Suppose that $\lambda$ is a multipartition of $n$. Then
$$\gamma_\tllam=q^{\len(w_\lambda)}\!\!\!\!\!\!\!%
\prod_{(i,j,s)\in[\lambda]}\!\!%
    \frac{[h^{\lambda^{(s)}}_{ij}]_q}{[\ell^{\lambda^{(s)}}_{ij}]_q}
    \prod_{t=s+1}^r\!\! (q^{j-i}Q_s-q^{\lambda^{(t)}_1}\!\!Q_t)
              \prod_{k=1}^{\lambda^{(t)}_1}
	            \frac{(q^{j-i}Q_s-q^{k-1-\lambda^{(t)'}_k}\!\!Q_t)}
	                 {(q^{j-i}Q_s-q^{k-\lambda^{(t)'}_k}\!\!Q_t)}.
$$
$($Note that $\lambda^{(t)'}_k$ is the length of the $k^\th$ column of
$[\lambda^{(t)}].)$
\label{gamma_tllam}
\end{Lemma}

\begin{proof}
We argue by induction on $n$. If $n=0$, both sides are~$1$
and there is nothing to prove (by convention, empty products are~$1$).
Suppose that $n>0$. Let $\mu=\Shape\(\tllam\rest(n-1)\)$; then
$\mu$ is a multipartition of $n-1$. Recall that
$\alpha(\lambda)=\frac12\sum_{s=1}^r
                    \sum_{i\ge1}(\lambda^{(s)}_i-1)\lambda^{(s)}_i$. 
Applying the definitions (see (\ref{gamma definition})),
$$\frac{\gamma_\tllam}{\gamma_{\t_\mu}}
  =q^{\len(w_\lambda)+\alpha(\lambda)-\len(w_\mu)-\alpha(\mu)}
  \frac{\prod_{x\in \mathscr A_\tllam(n)}\(\res_\tllam(n)-\res(x)\)}%
     {\prod_{y\in \mathscr R_\tllam(n)}\(\res_\tllam(n)-\res(y)\)}.$$
Assume that $n$ appears in row $a$ and column $b$ of
$\t_\lambda^{(c)}$.  First consider the contribution that the
addable and removable nodes in $[\lambda^{(c)}]$ make to
$\gamma_\tllam$. Looking at the definitions above 
(\ref{gamma definition}), these nodes occur in pairs $(x,y)$
where~$y<(a,b,c)$ is a removable node in row~$i$ and~$x$ is an addable
node in row~$i+1$ for some $i\ge a$. If $x$ is in column $d$ of
$[\lambda^{(c)}]$ and~$y$ is in column $d'$ then $d\le d'<b$ and
\begin{align*}
\frac{\res_\tllam(n)-\res_\tllam(x)}{\res_\tllam(n)-\res_\tllam(y)}
  &=\frac{q^{b-a}Q_c-q^{d-(i+1)}Q_c}{q^{b-a}Q_c-q^{d'-i}Q_c}
   =\prod_{j=d}^{d'}\frac{q^{b-a}-q^{j-(i+1)}}{q^{b-a}-q^{j-i}}\\
  &=\prod_{j=d}^{d'}\frac{q^{j-i-1}(q^{b-a-j+i+1}-1)}%
                         {q^{j-i}(q^{b-a-j+i}-1)}
   =\prod_{j=d}^{d'}\frac{q^{-1}(q^{h^{\lambda^{(c)}}_{aj}}-1)}
                       {(q^{h^{\lambda^{(c)}}_{aj}-1}-1)}\\
  &=\prod_{j=d}^{d'}\frac{q^{-1}[h^{\lambda^{(c)}}_{aj}]_q}
                       {[h^{\lambda^{(c)}}_{aj}-1]_q}.
\end{align*}
Therefore,
$$
\frac{\prod_{(i,j,c)\in \mathscr A_\tllam(n)}
                \(\res_\tllam(n)-\res_\tllam(i,j,c)\)}%
     {\prod_{(i,j,c)\in \mathscr R_\tllam(n)}
                \(\res_\tllam(n)-\res_\tllam(i,j,c)\)}
     =\prod_{j=1}^{b-1}\frac{q^{-1}[h^{\lambda^{(c)}}_{aj}]_q}
                       {[h^{\lambda^{(c)}}_{aj}\!-1]_q}
    =q^{1-b}\prod_{j=1}^{b-1}
       \frac{[h^{\lambda^{(c)}}_{aj}]_q}{[h^{\mu^{(c)}}_{aj}]_q}.
$$
Note that $\alpha(\lambda)=\alpha(\mu)+b-1$; hence, by induction,
this accounts for the factor $q^{\len(w_\lambda)}$ in Lemma~\ref{gamma_tllam}. 

Now $\gamma_{\t_\mu}$ is known by induction and it contains as a
factor the left hand term in the product below. Further,
$\ell^{\mu^{(s)}}_{ij}=\ell^{\lambda^{(s)}}_{ij}$ and
$h^{\mu^{(s)}}_{ij}=h^{\lambda^{(s)}}_{ij}$ if $(i,s)\ne(a,c)$,
and $\ell^{\mu^{(c)}}_{aj}=\ell^{\lambda^{(c)}}_{aj}$ for 
$1\le j<b$,~so
\begin{align*}
\bigg(\prod_{(i,j,s)\in[\mu]}%
         \frac{[h^{\mu^{(s)}}_{ij}]_q}{[\ell^{\mu^{(s)}}_{ij}]_q}\bigg)%
	 \bigg(\prod_{j=1}^{b-1} 
	 \frac{[h^{\lambda^{(c)}}_{aj}]_q}{[h^{\mu^{(c)}}_{aj}]_q}\bigg)
    &=\bigg(\prod_{\substack{(i,j,s)\in[\lambda]\\ (i,s)\ne(a,c)}}%
      \frac{[h^{\lambda^{(s)}}_{ij}]_q}{[\ell^{\lambda^{(s)}}_{ij}]_q}
     \bigg)\bigg(%
     \prod_{j=1}^{b-1} 
       \frac{[h^{\lambda^{(c)}}_{aj}]_q}{[\ell^{\lambda^{(c)}}_{aj}]_q}
     \bigg)\\
   &=\prod_{(i,j,s)\in[\lambda]}%
      \frac{[h^{\lambda^{(s)}}_{ij}]_q}{[\ell^{\lambda^{(s)}}_{ij}]_q},
\end{align*}
since $[h^{\lambda^{(c)}}_{ab}]_q=1=[\ell^{\lambda^{(c)}}_{ab}]_q$. 
This accounts for the left hand factor in the expression
for~$\gamma_\tllam$ given in the statement of the Lemma.

Finally, consider the nodes in 
$\mathscr A_\tllam(n)$ and $\mathscr R_\tllam(n)$ which are in
component~$t$ for some $t>c$ (there are no such nodes for $t<c$).
Again, almost all of the addable and removable nodes in component~$t$
occur in pairs placed in consecutive rows; however, this time there is
also an additional addable node at the end of the first row
of~$\lambda^{(t)}$. As above, it is easier to insert extra factors
which cancel out and so take a product over all of the columns
of~$\lambda^{(t)}$. An argument similar to that above shows that the
nodes in $\mathscr A_\tllam(n)$ and $\mathscr R_\tllam(n)$ which do
not belong to component~$c$ contribute the factor 
$$\prod_{t=c+1}^r (q^{b-a}Q_c-q^{\lambda^{(t)}_1}Q_t)
  \prod_{k=1}^{\lambda^{(t)}_1}
      \frac{(q^{b-a}Q_c-q^{k-1-\lambda^{(t)'}_k}\!\!Q_t)}
           {(q^{b-a}Q_c-q^{k-\lambda^{(t)'}_k}\!\!Q_t)}$$ 
to $\gamma_\tllam$. Using induction to combine the formulae above 
proves the Lemma.
\end{proof}

We can now give a closed formula for the Schur elements.

\begin{Corollary}
Suppose that $\lambda$ is a multipartition of $n$ and for
$1\le s<t\le r$ let
$$X^\lambda_{st}=\!\!\!\prod_{(i,j)\in[\lambda^{(t)}]}(q^{j-i}Q_t-Q_s)
\ \cdot\!\!\!\!
  \prod_{(i,j)\in[\lambda^{(s)}]}(q^{j-i}Q_s-q^{\lambda^{(t)}_1})
      \prod_{k=1}^{\lambda^{(t)}_1}
        \frac{(q^{j-i}Q_s-q^{k-1-\lambda^{(t)'}_k}Q_t)}
             {(q^{j-i}Q_s-q^{k-\lambda^{(t)'}_k}Q_t)}.$$
Then
$$s_\lambda(\Q)=(-1)^{n(r-1)}(Q_1\dots Q_r)^{-n}q^{-\alpha(\lambda')}
  \prod_{s=1}^r\prod_{(i,j)=[\lambda^{(s)}]}Q_s[h_{ij}^{\lambda^{(s)}}]_q
  \cdot\prod_{1\le s<t\le r}X^\lambda_{st}.$$
\label{closed}
\end{Corollary}

\begin{proof}
Looking at the definitions,
$[k]_q'=1+q^{-1}+\dots+q^{1-k}=q^{1-k}[k]_q$; therefore,
$([k]_q^!)'=q^{-\frac12k(k-1)}[k]_q^!$; consequently,
$([\lamp]_q^!)'=q^{-\alpha(\lamp)}[\lamp]_q^!$. Next
observe that
$[\lamp]_q^!
        =\prod_{(i,j,s)\in[\lambda]}[\ell^{\lambda^{(s)}}_{ij}]_q$.
Further, $(i,j)\in[\lambda^{(s)}]$ if and only if
$(j,i)\in[\lamp^{(r-s+1)}]$.  Therefore, applying ${}'$ to
(\ref{gamma properties})(i) (and swapping the roles of $s$ and $t$ in
the right hand factor),
$$\gamma_{\tlamp}'=q^{-\alpha(\lamp)}
 \bigg(\prod_{(i,j,s)\in[\lambda]}[\ell^{\lambda^{(s)}}_{ij}]_q\bigg)
 \bigg(\prod_{1\le s<t\le r}\prod_{(i,j)\in[\lambda^{(t)}]}
	            (q^{j-i}Q_t-Q_s)\bigg).$$
By Corollary~\ref{Schur elements},
$s_\lambda(\Q)=(-1)^{n(r-1)}q^{-\len(w_\lambda)}
   \gamma_\tllam\gamma_{\tlamp}'\prod_{s=1}^r Q_s^{|\lambda^{(s)}|-n}$; 
so the result now follows by Lemma~\ref{gamma_tllam}.
\end{proof}

\vspace{-5mm}
\begin{Example}
It is straightforward to check that Corollary~\ref{closed} gives the
same rational functions for the Schur elements $s_{\eta_t}(\Q)$ as 
were obtained in Example~\ref{n_t example}.
\end{Example}

It is not at all obvious that the formula in Corollary~\ref{closed} for the
Schur elements agrees with that conjectured by Malle, so we now show
that this is the case. This takes quite a lot of work; however, it
also results in a more symmetrical formula for $s_\lambda(\Q)$. To do
this we need to rewrite the Schur elements as functions of beta
numbers (or first column hook lengths). 

Define the length of a partition $\sigma$ to be the smallest integer
$\len(\sigma)$  such that $\sigma_i=0$ for all~$i>\len(\sigma)$.
The {\sf length} of a multipartition $\lambda$ is 
$\len(\lambda)=\max\set{\len(\lambda^{(s)})|1\le s\le r}$.

Fix an integer $L$ such that $L\ge\len(\lambda)$. The 
$L$--{\sf beta numbers} for~$\lambda^{(s)}$ are the integers
$\beta^{(s)}_i=\lambda^{(s)}_i+L-i$ for $i=1,\dots,L$; note that
$\beta^{(s)}_1>\dots>\beta^{(s)}_L\ge0$. Malle calls the $r\times L$
matrix $B=\(\beta^{(s)}_i)_{s,i}$ the $L$--{\sf symbol} of $\lambda$.
Actually, this is not quite one of Malle's symbols; it is what Brou\'e
and Kim~\cite{BK} call~$B$ an {\it ordinary symbol}. We also let
$B_s=\{\beta^{(s)}_1,\dots,\beta^{(s)}_L\}$ for $s=1,\dots,r$.

If we change $L$ to $L+1$ the beta set $B_s$ is {\it shifted} to
$\{\beta^{(s)}_1+1,\dots,\beta^{(s)}_L+1,0\}$. A function
of beta numbers is {\sf invariant under beta shifts} if it is unchanged by
such transformations; equivalently, the function is independent of 
$L$ provided that $L$ is large enough. For example, the formula for
$s_\lambda(\Q)$ below is invariant under beta shifts since
$s_\lambda(\Q)$ does not depend on~$L$.

\begin{Theorem}
Suppose that $\lambda$ is a multipartition of $n$ with
$L$--symbol $B=(\beta^{(s)}_i)_{s,i}$ such that $L\ge\len(\lambda)$. 
Then
$$s_\lambda(\Q)=(-1)^{a_{rL}}q^{b_{rL}}
\frac{\Prod_{1\le s<t\le r}(Q_s-Q_t)^L\ \cdot
      \prod_{1\le s,t\le r}\prod_{\alpha_s\in B_s}
      \prod_{1\le k\le\alpha_s} (q^kQ_s-Q_t)}
     {(q-1)^n(Q_1\dots Q_r)^n
      \Prod_{1\le s\le t\le r}
      \prod_{\substack{(\alpha_s,\alpha_t)\in B_s\times B_t\\
                      \alpha_s>\alpha_t\If s=t}}
             (q^{\alpha_s}Q_s-q^{\alpha_t}Q_t)},$$
where $a_{rL}=n(r-1)+\binom r2\binom L2$ 
and $b_{rL}=\frac{rL(L-1)(2rL-r-3)}{12}$.
\label{Malle}
\end{Theorem}

\begin{proof}
Adopting the notation of Corollary~\ref{closed},
$$s_\lambda(\Q)=(-1)^{n(r-1)}(Q_1\dots Q_r)^{-n}q^{-\alpha(\lambda')}
  \prod_{s=1}^r\prod_{(i,j)=[\lambda^{(s)}]}Q_s[h_{ij}^{\lambda^{(s)}}]_q
  \cdot\prod_{1\le s<t\le r}X^\lambda_{st}.$$
We consider each of these factors separately.

First, let $\sigma$ be a partition of $m$ with beta numbers
$(\beta_1,\dots,\beta_L)$, where $L\ge\len(\sigma)$. Then
$$\prod_{(i,j)\in[\sigma]}[h^\sigma_{ij}]_q
      =(q-1)^{-m}\prod_{(i,j)\in[\sigma]}(q^{h^\sigma_{ij}}-1)
      =(q-1)^{-m}\frac{\Prod_{i=1}^L\prod_{k=1}^{\beta_i}(q^k-1)}
                {\Prod_{1\le i<j\le L}(q^{\beta_i-\beta_j}-1)}.
$$
This is quite well--known and  is easily proved by first observing
that the right hand side is invariant under beta shifts and then by
arguing by induction on the number of columns.

Ignoring the leading term $(q-1)^m$ on the right hand side, 
the number of factors in the numerator is $m+\binom L2$, whereas
the number of factors in the denominator is $\binom L2$. Consequently,
multiplying the left hand side by $Q^m$, say, is the same as
multiplying each factor, top and bottom, on the right hand side by
$Q$.  Therefore, returning to the multipartition $\lambda$, 
$$q^{-\alpha(\lamp)}\prod_{s=1}^r\prod_{(i,j)\in[\lambda^{(s)}]}
            Q_s[h_{ij}^{\lambda^{(s)}}]_q
  =q^{N-\alpha(\lamp)}(q-1)^{-n}\prod_{s=1}^r\frac%
     {\Prod_{i=1}^L\prod_{k=1}^{\beta_i^{(s)}}(q^kQ_s-Q_s)}
     {\Prod_{1\le i<j\le L}(q^{\beta^{(s)}_i}Q_s-q^{\beta^{(s)}_j}Q_s)},
$$
where $N=\sum_{s=1}^r\sum_{j=1}^L(j-1)\beta^{(s)}_j$. It is
well--known (in the right circles) and easy enough to check that
$\alpha(\lamp)=\sum_{s=1}^r\sum_{j=1}^L(j-1)\lambda^{(s)}_j$;
therefore, $N-\alpha(\lamp)=r\binom L3$.

Fix $s$ and $t$ with $1\le s<t\le r$; there are $\binom r2$ such
choices. A quick calculation shows that $b_{rL}=r\binom L3+\binom
r2\sigma(L)$, where $\sigma(L)=(2L-1)L(L-1)/6=\sum_{j=1}^{L-1}j^2$.
Therefore, in order to complete 
the proof it is enough to show that $X^\lambda_{st}=Y^\lambda_{st}$ where 
$$Y^\lambda_{\s\t}=(-1)^{\binom L2}q^{\sigma(L)}
   \frac{(Q_s-Q_t)^L\!\!\!
     \Prod_{\alpha_s\in B_s}\prod_{1\le k\le\alpha_s}(q^kQ_s-Q_t)
         \ \cdot
         \prod_{\alpha_t\in B_t}\prod_{1\le k\le\alpha_t}(q^kQ_t-Q_s)}
        {\Prod_{(\alpha_s,\alpha_t)\in B_s\times B_t}
                     (q^{\alpha_s}Q_s-q^{\alpha_t}Q_t)},$$
It is not hard to see that $Y^\lambda_{st}$ is invariant under beta
shifts so we may change $L$ arbitrarily, provided that
$L\ge\max\big\{\len(\lambda^{(s)}), \len(\lambda^{(t)})\big\}$.  We
prove our claim that $X^\lambda_{st}=Y^\lambda_{st}$ by induction in
three incremental steps. We start the induction by observing that both
products are equal to~$1$ when $\lambda^{(s)}=\lambda^{(t)}=(0)$ ---
to see this it is easiest to take $L=0$.

Next consider the case where $\lambda^{(t)}=(0)$. Since
$Y^\lambda_{st}$ is invariant under beta shifts we may assume that
$\len(\lambda^{(s)})=L$. Assume by way of induction that we have
proved the claim for $(\lambda^{(s)}_1,\dots,\lambda^{(s)}_{L-1})$.
Adding a non--empty $L^{\text{th}}$ row to $\lambda^{(s)}$ changes
$X^\lambda_{st}$ by the factor
\begin{align*}
\prod_{k=1-L}^{\lambda^{(s)}_L-L}(q^kQ_s-Q_t)
    &=\prod_{k=1}^{L-1}q^{-k}(Q_s-q^kQ_t)
      \ \cdot(Q_s-Q_t)\cdot\!\!
      \prod_{k=1}^{\lambda^{(s)}_L-L}(q^kQ_s-Q_t)\\
    &=(-1)^{L-1}
      \frac{
        \Prod_{1\le k\le L-1}\!\!\!\!\!(q^kQ_t-Q_s)
        \cdot(Q_s-Q_t)\cdot\!\!\!\!\!\!\!
        \prod_{1\le k\le \lambda^{(s)}_L}\!\!\!\!\!(q^kQ_s-Q_t)}
           {q^{\binom L2}
            \Prod_{\lambda^{(s)}_L-L+1\le k\le\lambda^{(s)}_L}
                       (q^kQ_s-Q_t)}\\
    &=(-1)^{L-1}
      \frac{
        \Prod_{1\le k\le\beta^{(t)}_L}\!\!\!\!\!(q^kQ_t-Q_s)
        \cdot(Q_s-Q_t)\cdot\!\!\!\!\!\!\!
        \prod_{1\le k\le \beta^{(s)}_L}\!\!\!\!\!(q^kQ_s-Q_t)}
         {\Prod_{\alpha_t\in B_t} 
              (q^{\beta^{(s)}_L}Q_s-q^{\alpha_t}Q_t)},
\end{align*}
the last equality following because $\beta^{(s)}_L=\lambda^{(s)}_L$
and $\beta^{(t)}_i=L-i$ for $i=1,\dots,L$. Notice that the
$(L-1)$--beta numbers for
$(\lambda^{(s)}_1,\dots,\lambda^{(s)}_{L-1})$ all increase by~$1$ when
we add the extra row~$\lambda^{(s)}_L$ to $\lambda^{(s)}$. Let $B_s'$
be the set of $(L-1)$--beta numbers
for~$(\lambda^{(s)}_1,\dots,\lambda^{(s)}_{L-1})$ and 
let~$B_t'=\{0,1,\dots,L-2\}$ be the set of $(L-1)$--beta numbers
for~$\lambda^{(t)}$. Then
$$\begin{array}[c]{ll}
\MC2l{\Prod_{(\alpha_s,\alpha_t)\in B_s\times B_t}
   (q^{\alpha_s}Q_s-Q^{\alpha_t}Q_t)
   =\prod_{1\le i\le L}\prod_{0\le k<L}
       (q^{\beta^{(s)}_i}Q_s-q^kQ_t)}\\
\qquad&=\Prod_{\alpha_t\in B_t}(q^{\beta^{(s)}_L}Q_s-q^{\alpha_t}Q_t)
  \cdot\!\!\!\prod_{1\le i<L}(q^{\beta^{(s)}_i}Q_s-Q_t)
  \cdot\!\!\!\!\!\!\!\!\!\!
  \prod_{(\alpha_s,\alpha_t)\in B_s'\times B_t'}\!\!\!\!
     q(q^{\alpha_s}Q_s-q^{\alpha_t}Q_t).
\end{array}$$
This equation allows us to rewrite the denominator of the preceding
equation and so see that the change in $X^\lambda_{st}$ is the same as the
change in $Y^\lambda_{st}$ (in particular, the change of the scalar is
$(-1)^{L-1}q^{(L-1)^2}$ in both cases). This proves our claim when
$\lambda^{(t)}=(0)$.

The next step is to fix $\lambda^{(s)}$ and assume that
$\lambda^{(t)}=(a)$ for some $a\ge0$. If $\lambda^{(s)}=(0)$ in this
case then it is straightforward to check the claim (or to modify the
argument below), so assume that $\lambda^{(s)}\ne(0)$ and let
$L=\len(\lambda^{(s)})\ge1\ge\len(\lambda^{(t)})$. The case $a=0$ we
already understand.  Next, changing $\lambda^{(t)}$ from $(a-1)$ to
$(a)$ changes $X^\lambda_{st}$ by the factor
$$\begin{array}{ll}
\MC2l{(q^{a-1}Q_t-Q_s)\Prod_{(i,j)\in[\lambda^{(s)}]}
          \frac{(q^{j-i}Q_s-q^aQ_t)}{(q^{j-i}Q_s-q^{a-1}Q_t)}
          \frac{(q^{j-i}Q_s-q^{a-2}Q_t)}{(q^{j-i}Q_s-q^{a-1}Q_t)}}\\
\quad&=(q^{a-1}Q_t-Q_s)\Prod_{(i,j)\in[\lambda^{(s)}]}
          \frac{(q^{a-1}Q_t-q^{j-i-1}Q_s)(q^{j-i+1}Q_s-q^{a-1}Q_t)}
               {(q^{j-i}Q_s-q^{a-1}Q_t)(q^{a-1}Q_t-q^{j-i}Q_s)}\\
     &=(q^{a-1}Q_t-Q_s)\Prod_{i=1}^L
 \frac{(q^{a-1}Q_t-q^{-i}Q_s)}{(q^{a-1}Q_t-q^{1-i}Q_s)}
 \frac{(q^{\lambda^{(s)}_i-i+1}Q_s-q^{a-1}Q_t)}
      {(q^{\lambda^{(s)}_i-i}Q_s-q^{a-1}Q_t)}\\
 &=(q^{a-1}Q_t-q^{-L}Q_s)\Prod_{i=1}^L
    \frac{(q^{\beta^{(s)}_i+1}Q_s-q^{\beta^{(t)}_1}Q_t)}
         {(q^{\beta^{(s)}_i}Q_s-q^{\beta^{(t)}_1}Q_t)}\\
 &=(q^{\beta^{(t)}_1}Q_t-Q_s)\Prod_{i=1}^L
    \frac{(q^{\beta^{(s)}_i}Q_s-q^{\beta^{(t)}_1-1}Q_t)}
         {(q^{\beta^{(s)}_i}Q_s-q^{\beta^{(t)}_1}Q_t)}.\\
\end{array}$$
This last product is exactly the change in $Y^\lambda_{st}$ so 
we now know that $X^\lambda_{st}=Y^\lambda_{st}$ when $\lambda^{(t)}=(a)$.

Finally, suppose that $\lambda^{(t)}$ has more than one row. For
convenience, let $l=\len(\lambda^{(s)})$, $m=\len(\lambda^{(t)})>1$
and $b=\lambda^{(t)}_m$ and assume that $L\ge l,m$. When we add row
$m$ to $\lambda^{(t)}$, $X^\lambda_{st}$ changes by the factor
$$\begin{array}{ll}
\MC2l{\Prod_{j=1}^b(q^{j-m}Q_t-Q_s)\ \cdot\!\!\!
  \prod_{(i,j)\in[\lambda^{(s)}]}\prod_{k=1}^b
     \frac{(q^{j-i}Q_s-q^{k-m-1}Q_t)}{(q^{j-i}Q_s-q^{k-m}Q_t)}
     \frac{(q^{j-i}Q_s-q^{k-m+1}Q_t)}{(q^{j-i}Q_s-q^{k-m}Q_t)}}\\
\  &=\Prod_{j=1}^b(q^{j-m}Q_t-Q_s)\cdot\!\!\!\!\!\!
  \prod_{(i,j)\in[\lambda^{(s)}]}\prod_{k=1}^b
     \frac{(q^{j-i+1}Q_s-q^{k-m}Q_t)}{(q^{j-i}Q_s-q^{k-m}Q_t)}
     \frac{(q^{k-m}Q_t-q^{j-i-1}Q_s)}{(q^{k-m}Q_t-q^{j-i}Q_s)}\\
     &=\Prod_{j=1}^b(q^{j-m}Q_t-Q_s)\cdot
  \prod_{i=1}^l\prod_{k=1}^b
  \frac{q(q^{\lambda^{(s)}_i-i}Q_s-q^{k-m-1}Q_t)}
       {(q^{\lambda^{(s)}_i-i}Q_s-q^{k-m}Q_t)}
  \frac{(q^{k-m}Q_t-q^{-i}Q_s)}{(q^{k-m}Q_t-q^{1-i}Q_s)}\\
     &=\Prod_{j=1}^b(q^{j-m}Q_t-Q_s)\cdot
q^{lb}
\prod_{i=1}^l\frac{(q^{\lambda^{(s)}_i-i}Q_s-q^{-m}Q_t)}
                  {(q^{\lambda^{(s)}_i-i}Q_s-q^{b-m}Q_t)}
\cdot
\prod_{k=1}^b\frac{(q^{k-m}Q_t-q^{-l}Q_s)}{(q^{k-m}Q_t-Q_s)}\\
     &=\Prod_{k=1}^b(q^{k-m+l}Q_t-Q_s)\cdot
\prod_{i=1}^l\frac{(q^{\beta^{(s)}_i}Q_s-q^{\beta^{(t)}_m-b}Q_t)}
                  {(q^{\beta^{(s)}_i}Q_s-q^{\beta^{(t)}_m}Q_t)}\\
     &=\Prod_{k=L-m+1}^{\beta^{(t)}_m}(q^kQ_t-Q_s)\cdot
\prod_{i=1}^L\frac{(q^{\beta^{(s)}_i}Q_s-q^{\beta^{(t)}_m-b}Q_t)}
                  {(q^{\beta^{(s)}_i}Q_s-q^{\beta^{(t)}_m}Q_t)},
\end{array}$$
where the last line follows by a short calculation since
$\beta^{(t)}_m=b+L-m$ (and $l\le L$). As $\beta^{(t)}_m$ is the only
beta number that has changed, this factor is precisely the change in
$Y^\lambda_{st}$ when an extra row is added to $\lambda^{(t)}$. 

We have now shown that $X^\lambda_{st}=Y^\lambda_{st}$ in all cases,
so the theorem is proved.  
\end{proof}

Finally, comparing Theorem~\ref{Malle} with \cite[Prop 3.17]{BK} we
see that our formula for the Schur elements agrees with Malle's --- as
it must because Malle's conjecture has already been proved by Geck,
Iancu and Malle~\cite{GIM}. Actually, there is still a small amount of
work to be done in reconciling the two formulas because Brou\'e and
Kim~\cite{BK} write the exponent of $q$ as a sum of binomial
coefficients; however, their expression simplifies to give $b_{rL}$.

\section*{\it Acknowledgements}

\noindent
I would like to thank the referee for many useful comments and
corrections and Jean Michel for writing some {\sc Chevie} code for
computing the Schur elements.

\let\em\it


\end{document}